\documentclass[twoside,12pt]{article}
\usepackage{amsmath}
\usepackage{amssymb}
\usepackage{amscd}
\usepackage{amsthm}%+

\numberwithin{equation}{section}

\theoremstyle{plain}
\newtheorem{theorem}{Theorem}[section]

\newtheorem{Le}[theorem]{Lemma}

\theoremstyle{remark}

\theoremstyle{definition}

\newcommand{\R}{\mathbb R}
\newcommand{\C}{\mathbb C}

\newcommand{\Q}{\mathbb Q}

\def\ra{\rightarrow}

\def\e{\emph}
\def\i{\infty}

\def\b{\begin}

\begin{document}

\title{    %\flushleft
{\bf{Quasiconformal maps on model Filiform groups}} }
\author{Xiangdong Xie\footnote{Partially supported by NSF grant DMS--1265735.}}
  %\author{Xiangdong Xie}
\date{  }
\maketitle

\begin{abstract}
We   describe all quasiconformal maps on   the higher  (real and complex)  
     model Filiform  groups  equipped with the Carnot metric,  including  non-smooth ones.  
  %  generalizing a result of Ben Warhurst \cite{W}  to the nonsmooth setting.    
 These maps  have very special forms. 
  In particular, they are all biLipschitz  and preserve multiple foliations. 
  The results in this paper have implications to the large scale geometry of 
     nilpotent   Lie groups  and  
negatively curved solvable Lie groups.

\end{abstract}

%\noindent
%{\bf{Mathematics Subject Classification(2000).}} Primary 53C23, 51F99, 57M20;
% Secondary  53C70, 57M60.
% \newline
%{\bf{Keywords.}}Key words:  quasisymmetric  homeomorphism, negative curvature,
 %solvable Lie group, visual metric,parabolic visual metric.

%Tits boundary, Tits metric, CAT(0), 2-complex, geodesic.

%\keywords}    test  \endkeywords
%\subjclass  51  \endsubjclass

%\tableofcontents

% negatively curved solvable Lie groups.

%\keywords    test  \endkeywords
%\subjclass  51  \endsubjclass

%\tableofcontents

%\vspace{3mm} \noindent
 {\small {\bf{Mathematics Subject
Classification (2010).}}  22E25,   30L10,   53C17.

  %53C24,  53C17,  53C23, 53C30

%20F65,  30C65, 53C20.

%\noindent
%{\bf{Mathematics Subject Classification(2000).}} Primary 53C23, 51F99, 57M20;
% Secondary  53C70, 57M60.
% \newline
%{\bf{Keywords.}} Tits boundary, Tits metric, CAT(0), 2-complex, geodesic.

% 	22E25   	Nilpotent and solvable Lie groups

%\vspace{3mm} \noindent {\small {\bf{Mathematics Subject
%Classification (2000).}} 54F45, 30C65.
%20F65. %57M20, 20F67, 20E07.}

%53C24       Rigidity results
%53C17       Sub-Riemannian geometry
%53C23       Global geometric and topological methods (�� la Gromov);
 %differential geometric analysis on metric spaces
%53C30       Homogeneous manifolds
%30L10       Quasiconformal mappings in metric spaces

%30C65  (1991-now) Quasiconformal mappings in ${\bf R}^n$, other generalizations
%20F34 Fundamental groups and their automorphisms
%20F65 Geometric group theory
%20F67 Hyperbolic groups and nonpositively curved groups
%20F99 None of the above, but in this section
%20E07 Subgroup theorems; subgroup growth
%53C20 Global Riemannian geometry, including pinching
%53C23 Global topological methods
%54F45 Dimension theory
%57M20 Two-dimensional complexes
%57M60 Group actions in low dimensions
%(1) (2) (a) (b) (i) (ii)
%20F69 Asymptotic properties of groups

%$\{x_i\}_{i=1}^\i$  converges to $\xi\in \ol{X}$
% an index two subgroup

               %\vspace{3mm} \noindent {\small {\bf{Key words.}} }

%hyperbolic element, parabolic element,
%quasi-convex.}

\setcounter{section}{0} \setcounter{subsection}{0}

\section{Introduction}\label{s0}

In this paper we study quasiconformal maps on the higher   real and complex model Filiform groups
  equipped with the Carnot metric.
  We   identify
     all such maps.  They are all biLipschitz  and preserve multiple foliations.  
    We do not impose any regularity conditions on the  quasiconformal maps.
      However,     the group structure forces rigidity and regularity.  In particular, in the case of   higher complex 
  model   Filiform groups,   up to taking complex conjugation  all quasiconformal maps are biholomorphic.

Let $K$ be a field. We only consider the case when $K$ is $\R$ or $\C$.  
 The $n$-step ($n\ge 2$)  model Filiform algebra $\mathcal{F}_K^n$   over $K$ is an $(n+1)$-dimensional  Lie 
 algebra over $K$.  
  It has a   basis $\{e_1, e_2, \cdots, e_{n+1}\}$ 
  and  the only non-trivial bracket relations are 
    $[e_1, e_j]=e_{j+1}$  for $2\le j\le n$.    
    The Lie algebra   $\mathcal{F}_K^n$   admits a  direct sum decomposition of vector subspaces
  $\mathcal{F}_K^n=V_1\oplus \cdots \oplus V_{n}$, where $V_1$ is the  linear subspace spanned by $e_1, e_2$, and $V_j$ ($2\le j\le n$)  is the linear subspace spanned by $e_{j+1}$.     It is easy to  check that 
  $[V_1, V_j]=V_{j+1}$ for $1\le j\le n$,  where $V_{n+1}=\{0\}$.   Hence $\mathcal{F}_K^n$  is a 
  graded      Lie algebra. 
 For $K=\R$ or $\C$,    the connected and simply connected Lie group with Lie algebra 
$\mathcal{F}_K^n$    will be denoted by  $F_K^n$  and is called the 
   $n$-step   model Filiform  group over $K$.

The two dimensional subspace $V_1$  of 
$\mathcal{F}_\R^n$    determines a left invariant distribution (so called horizontal distribution)  on  $F_\R^n$.  On $V_1$, we consider the inner product with $e_1$, $e_2$ as orthonormal basis. 
 This norm on $V_1$ then induces a Carnot metric  $d_c$ on    $F_\R^n$.
  Similarly   the first layer $V_1$ of  $\mathcal{F}_\C^n$   is a 
 four dimensional real vector subspace spanned by $e_1, ie_1, e_2, ie_2$ ($i=\sqrt{-1}$),  and   it 
 determines a left invariant distribution   on  $F_\C^n$.  On $V_1$ of   $\mathcal{F}_\C^n$, we consider the inner product with $e_1, ie_1$, $e_2, ie_2$ as orthonormal basis. 
 This norm on $V_1$ then induces a Carnot metric  $d_c$ on    $F_\C^n$.

Recall that, for a simply connected nilpotent Lie group  $G$ with Lie algebra $\mathcal G$, the exponential map
   $\text{exp}:  \mathcal{G}\ra G$ is  a  diffeomorphism.    We shall identify $\mathcal{G}$ and $G$ via the exponential map and denote the group operation  by $*$.    Notice that every element 
$x\in \mathcal{F}_\R^n$  can be uniquely written as $x=x_1e_1*(x_2 e_2+\cdots  + x_{n+1} e_{n+1})$,
 where $x_i\in \R$, $1\le i\le n+1$. 

Fix $n\ge 2$.  Let $h: \R\ra \R$ be a   Lipschitz   function. 
    Set $h_2=h$  and define $h_j:  \R  \ra \R$,  $3\le j\le n+1$, inductively as follows:
  $$h_j(x)=-\int_0^x h_{j-1}(s) ds.$$
    Define $F_h:   F_\R^n\ra    F_\R^n$   by:
    $$F_h\big(x_1e_1*\sum_{j=2}^{n+1}x_je_j\big)=x_1e_1*\sum_{j=2}^{n+1}\big(x_j+h_j(x_1)\big)e_j.$$

\b{Th}\label{thA}
  Let $n\ge 3$.  
    A  homeomorphism  $F: (F_\R^n, d_c)  \ra    (F_\R^n,  d_c)$   is   quasiconformal if and only if  it 
 is a finite composition of left translations, graded isomorphisms and maps of   the 
 form   $F_h$, where $h: \R\ra \R$ is a Lipschitz   function.

\end{Th}

  We remark that Ben Warhurst   \cite{W}    previously proved    a  similar statement under the assumption that the map $F$ is smooth.  We do not impose any regularity assumptions on the quasiconformal maps. 
    Theorem \ref{thA}  also provides   lots of examples of non-smooth quasiconformal  maps on 
  $F_\R^n$.

  The quasiconformal   maps on the complex Filiform groups are even more rigid:

\b{Th}\label{thB}
  Let $n\ge 3$.  
    A  homeomorphism  $F: (F_\C^n,  d_c)  \ra    (F_\C^n,  d_c)$  is a quasiconformal map if and only if it 
 is a finite composition of left translations  and  graded isomorphisms.

\end{Th}

   Model   Filiform groups are an important class of the so-called non-rigid   Carnot groups  \cite{W},  
    \cite{O}, \cite{OW}.      
    Recall that a  Carnot group is rigid if the space of smooth contact maps is finite dimensional, and called non-rigid otherwise.   % In fact, every non-rigid Carnot group contains a copy of the real or complex model Filiform group
   %  \cite{O}, \cite{OW}.    
Our results show that   even on these non-rigid Carnot groups  
quasiconformal maps  are  rigid in the sense that  they  are biLipschitz and  have very special forms.  

 When $n=2$, the group $F_\R^2$ is simply the first Heisenberg group. 
  It is well-known that quasiconformal maps on the   Heisenberg groups are very flexible.
 For instance, there   exist  quasiconformal maps between Heisenberg groups that change the Hausdorff dimension of certain subsets, see  \cite{B}.     There also exist biLipschitz maps of the Heisenberg groups that map   vertical lines to non-vertical curves \cite{X1}.

The group $F_\C^2$ is the first complex Heisenberg group. 
Recall that the $n$-th  complex Heisenberg algebra
  $\mathcal{H}^n_\C$  
 is  an    $(2n+1)$-dimensional complex Lie algebra and   has a complex vector space basis
  $X_i, Y_i, Z$ ($1\le  i\le n$)  with   the only non-trivial bracket relations   
  $[X_i, Y_i]=Z$,   $1\le i\le n$.  The Lie algebra $\mathcal{H}^n_\C$  is a 2-step Carnot algebra.  
The first layer $V_1$ of 
  $\mathcal{H}^n_\C$   is spanned by the $X_i, Y_i$, $1\le i\le n$ and has   complex  dimension $2n$. The second layer  $V_2$   is spanned by $Z$ and has complex  dimension $1$.   The $n$-th 
complex Heisenberg group $H^n_\C$ 
 is the simply connected  nilpotent   Lie group whose Lie algebra is    $\mathcal{H}^n_\C$.
  We always equip  $H^n_\C$  with a Carnot metric  associated with  $V_1$.  
We identify   both  $H^n_\C$   and  its  Lie algebra   $\mathcal{H}^n_\C$   with $\C^{2n+1}$.  
  So $V_1=\C^{2n}\times \{0\}$ and $V_2=\{0\}\times \C$.

  %Of course, it has  more bracket relations as a real Lie albegra coming from the fact that the bracket is complex linear:   $[X_j, iY_j]=[iX_j,  y_j]=i Z$  and $[iX_j,  iY_j]=-Z$.
  %one has bracket relations like $[X_i, iY_i]=i Z$ that come from the fact the bracket is complex linear.  
 %The first layer $V_1$ of 
 % $%\mathcal{H}^n_\C$   is spanned by the $X_i, Y_i$, $1\le i\le n$ and has   complex  dimension $2n$. The second layer  $V_2$   is spanned by $Z$ and has complex  dimension $1$. 
    %We identify $H^n_\C$   and  its  Lie algebra   $\mathcal{H}^n_\C$   with $\C^{2n+1}$.  
  %So $V_1=\C^{2n}\times \{0\}$ and $V_2=\{0\}\times \C$.    Let $\pi_1:  
%{H}^n_\C=\C^{2n+1}   \ra \C^{2n}$ be the projection onto $V_1$.  
  %A map $F: \C^{2n+1}   \ra \C^{2n+1} $  is called a lifting of a map $f: \C^{2n }  \ra \C^{2n}$
  %if $F(\pi_1^{-1}(p))=\pi_1^{-1}(f(p))$  for all $p\in \C^{2n}$.   
 It seems   complex Heisenberg groups are   also very rigid with respect to quasiconformal maps. In fact,
  using  the result of     Reimann  and   Ricci  \cite{RR}  one can show that if a quasiconformal map 
 $F:  H_\C^n\ra H_\C^n$  is also a $C^2$ diffeomorphism, then it is complex affine after possibly taking
   complex  conjugation,  see   Proposition \ref{affine}.  % Here the metric 
%on   $F_\C^2$  is a Carnot metric and $F_\C^2$  is identified with $\C^3$.   

 % We have the following:

\b{conjecture}\label{conjectyre-ch}
A    homeomorphism  $F:  H^n_\C\ra H^n_\C$ of the $n$-th  complex Heisenberg group is a  quasiconformal map
 if and only if it is a finite composition of left translations  and  graded isomorphisms. 
\end{conjecture}

For general Carnot groups, we have the following:

\b{conjecture}\label{qcbilip}
  Let $G$ be a Carnot group equipped with a  Carnot metric. If 
 $G$ is not an Euclidean group or an Heisenberg group, then every
   quasiconformal map $F: G\ra G$  is biLipschitz.
\end{conjecture}

  The results in this paper have implications for the large scale geometry of
nilpotent  groups  and  
 negatively curved 
 homogeneous manifolds. % and nilpotent Lie groups. 
  Each Carnot group arises  as  the (one point complement  of) ideal boundary  of some negatively
  curved   homogeneous manifold  \cite{H}.  
   Our results imply that each quasiisometry of the negatively curved homogeneous manifold assocaited to the higher model Filiform group is a rough isometry, that is, it must preserve the distance up to an additive constant.       Furthermore,  each quasiisometry between finitely generated nilpotent groups
  descends to a biLipschitz map between the asymptotic cones, which are Carnot groups.  
    Our results  say  that these biLipschitz maps preserve multiple foliations.  So 
 they provide information about the structure of the quasiisometries, at least after passing to the asymptotic cones.

The ideas in this  paper  can be used to show that quasiconformal maps on many Carnot groups are biLipschitz.    These results will  appear in  a forthcoming paper \cite{X2}.

In Section \ref{s1}, we  recall the basics about Carnot groups and the definition of model Filiform groups.
 In Section \ref{s2} we  study quasiconformal maps on the real model Filiform groups and  prove Theorem \ref{thA}.  In Section \ref{s4} we  consider complex Heisenberg groups; in particular, we prove a special case of Conjecture   \ref{conjectyre-ch}      which will be used later in Section \ref{s3}. 
 In    Section \ref{s3} we prove the rigidity result about quasiconformal   maps on the higher complex model Filiform groups (Theorem \ref{thB}).

\noindent {\bf{Acknowledgment}}. {This work was completed while the author was attending the workshop 
\lq\lq Interactions between analysis and geometry" at IPAM,  University of California at Los Angeles  from March to June 2013.  I would like to thank IPAM  for   financial support, excellent working conditions and  
   conducive atmosphere.
% condusive  atmorsphere.  
 I also would like to thank David Freeman and Tullia Dymarz for discussions about Carnot groups.  }
%   The author  is  partially supported  by NSF grant  DMS-1265735.  }

\section{Preliminaries}\label{s1}

In this Section we collect definitions and results that shall be needed later. 
  We first recall the basic definitions related to Carnot groups  in Subsection \ref{basics}.
   Then we  review    the 
   definition of  model Filiform groups (Subsection \ref{fili}),  the BCH formula (Subsection \ref{BCH}),
 the definitions of  quasi-similarity and quasisymmetric maps  (Subsection \ref{maps}),  and  
  Pansu differentiability theorem (Subsection \ref{pansud}).
  %  and the definitions of  quasi-similarity and quasisymmetric maps  (Subsection \ref{maps}).

\subsection{The basics}\label{basics}

A   \e{Carnot Lie algebra} is a finite dimensional Lie algebra
$\mathcal G$   together with  a direct sum   decomposition  
    $\mathcal
G=V_1\oplus V_2\oplus\cdots \oplus V_r$ 
  of   non-trivial   vector subspaces
 such that $[V_1,
V_i]=V_{i+1}$ for all $1\le i\le r$,
    where we set $V_{r+1}=\{0\}$.  The integer $r$ is called the
    degree of nilpotency of $\mathcal
G$. Every Carnot algebra
 $\mathcal
G=V_1\oplus V_2\oplus\cdots \oplus V_r$  admits a one-parameter
 family of automorphisms $\lambda_t: \mathcal
G \ra \mathcal G$, $t\in (0, \i)$, where
 $\lambda_t(x)=t^i x$ for  $x\in V_i$.
  Let   $\mathcal
G=V_1\oplus V_2\oplus\cdots \oplus V_r$
    and $\mathcal
G'=V'_1\oplus V'_2\oplus\cdots \oplus V'_s$  be two  Carnot
    algebras.
  A Lie algebra  homomorphism
  $\phi: \mathcal
G\ra \mathcal G'$
     is graded if $\phi$ commutes with $\lambda_t$ for
  all $t>0$; that is, if $\phi\circ \lambda_t=\lambda_t\circ
  \phi$.  We observe that $\phi(V_i)\subset V'_i$ for all $1\le i\le
  r$.

A  simply connected nilpotent Lie group is a  \e{Carnot group}
 if its Lie algebra is a Carnot algebra.
   %A Carnot group is reducible if its Lie algebra is reducible.
    % Equivalently,  a  Carnot group is reducible if it is isomorphic
     %to the direct product of two  Carnot groups.
    Let $G$ be a Carnot group with Lie algebra
      $\mathcal G=V_1\oplus \cdots \oplus V_r$.  The subspace $V_1$ defines
      a left invariant distribution  $H G\subset TG$ on $G$.    We fix a left invariant inner product on
          $HG$.
           An
      absolutely continuous curve $\gamma$ in $G$  whose velocity vector
       $\gamma'(t)$  is contained in  $H_{\gamma(t)} G$ for a.e. $t$
        is called  a horizontal curve.
          By Chow's theorem ([BR, Theorem 2.4]),   any two points
  of $G$ can be  connected by horizontal curves. Let $p,q\in G$, the
  \e{Carnot   metric} $d_c(p,q)$   between them is defined as
  the infimum of length of horizontal curves that join $p$ and $q$.

  Since the inner product on   $HG$ is left invariant, the Carnot
  metric on $G$ is also left invariant.  Different choices of inner
  product on $HG$ result in Carnot metrics that are biLipchitz
  equivalent.
    The Hausdorff dimension of $G$ with respect to  a  Carnot metric
    is  given by $\sum_{i=1}^r i\cdot \dim(V_i)$.

Recall that, for a simply connected nilpotent Lie group $G$ with Lie
algebra $\mathcal G$, the exponential map
  $\exp: {\mathcal G}\ra G$ is a diffeomorphism.  Under this identification the Lesbegue   measure on 
  $\mathcal G$  is a  Haar measure on $G$.  
 Furthermore, the
  exponential map induces %a Lie subalgebras of $\mathcal G$  and
  a  one-to-one correspondence between
    Lie subalgebras of $\mathcal G$   and
  connected Lie subgroups of $G$.
    %and

 It is often   more  convenient to work with homogeneous distances defined using norms
 than with Carnot metrics.  Let 
$\mathcal
G=V_1\oplus V_2\oplus\cdots \oplus V_r$  be a  Carnot algebra. 
  Write $x\in \mathcal G$ as $x=x_1+\cdots+ x_r$ with $x_i\in V_i$.  
  Fix a norm   $|\cdot|$ on each layer.  Define a norm $||\cdot||$ on $\mathcal G$
   by:
$$||x||=\sum_{i=1}^r |x_i|^{\frac{1}{i}}.$$
 Now define a homogeneous distance on $G=\mathcal G$  by:  $d(g,h)=||(-g)*h||$. 
  An important fact is that $d$ and $d_c$ are biLipschitz equivalent.   That is, there is a constant $C\ge 1$ such that    $d(p,q)/C\le d_c(p,q)\le C\cdot   d(p,q)$ for all $p, q\in G$. 
  It is often possible to  calculate  or estimate  $d$  by using  the BCH formula (see  Subsection \ref{BCH}).   
     Since we are only concerned with quasiconformal maps and biLipschitz maps, 
       it does not matter whether we use $d$ or $d_c$.

Let $G$ be a Carnot group with Lie algebra
      $\mathcal G=V_1\oplus \cdots \oplus  V_r$.
         Since $\lambda_t: \mathcal G\ra  \mathcal G$ ($t>0$) is a Lie algebra
         automorphism  and $G$ is simply connected,  there is a
         unique  Lie group automorphism  $\Lambda_t: G\ra G$ whose
         differential  at the identity is $\lambda_t$.
         %For $r>0$, the dilation $\delta_r:  \mathcal G\ra \mathcal
        %G$ is defined by
         % $\delta_r(\sum_{j=1}^k v_j)=\sum_{j=1}^k  r^j v_j$, where
           %$v_j\in V_j$.
          % Since the exponential map
           %$\exp: \mathcal G\ra G$ is a diffeomorphism,  the map
            % $\Lambda_t:=\exp  \circ \lambda_t\circ \exp^{-1}:
             %G\ra G$ is well-defined.
  %Notice that  $\Lambda_t$  is an automorphism of $G$.
      %Furthermore,
       For each $t>0$,
      $\Lambda_t$ is a
             similarity with respect to the Carnot metric:  $d(\Lambda_t(p),
             \Lambda_t(q))=t\, d(p,q)$ for any two points $p, q\in
             G$.  A  Lie group homomorphism
                $f: G\ra G'$ between two Carnot groups is a
                  graded homomorphism if  it commutes with
                    $\Lambda_t$ for all  $t>0$; that is, if
                      $f\circ \Lambda_t=\Lambda_t\circ  f$.
                         Notice that,
   a Lie group homomorphism
                $f: G\ra G'$ between two Carnot groups
                 is graded if and only if the corresponding Lie
                 algebra homomorphism is graded.

\subsection{The model Filiform groups}\label{fili}

Let $K$ be a field. We only consider the case when $K$ is $\R$ or $\C$.  
 The $n$-step ($n\ge 2$)  model Filiform algebra $\mathcal{F}_K^n$   over $K$ is an $(n+1)$-dimensional  Lie 
 algebra over $K$.  
  It has a   basis $\{e_1, e_2, \cdots, e_{n+1}\}$ 
  and  the only non-trivial bracket relations are 
    $[e_1, e_j]=e_{j+1}$  for $2\le j\le n$.    
  When $n=2$,  $\mathcal{F}_K^2$   is simply   the Heisenberg algebra   over $K$.   
  When $K=\C$,   $\mathcal{F}_\C^n$  can also be viewed as a real Lie algebra.  
   Since the brackets in  $\mathcal{F}_\C^n$  are complex linear,
   $\mathcal{F}_\C^n$   (when viewed  as a  real Lie algebra) 
     has  the following   additional   non-trivial bracket relations:  
 $[ie_1, e_j]=[e_1, ie_j]=ie_{j+1}$,  
 $[ie_1, ie_j]=-e_{j+1}$  for $2\le j\le n$.

For $K=\R$ or $\C$, the $n$-step ($n\ge 2$)  model Filiform group  ${F}_K^n$   over $K$
  is the simply connected Lie group whose Lie algebra is $\mathcal{F}_K^n$.
  For  $F^n_\R$,   we use the Carnot metric   corresponding to the inner product on $V_1$ with $e_1$ and $e_2$ as orthonormal basis.   The homogeneous distance on $F^n_\R$ is determined by the following norm:
  $$||\sum_{i=1}^{n+1} x_i e_i||=(x_1^2+x_2^2)^{\frac{1}{2}} +\sum_{i=3}^{n+1}|x_i|^{\frac{1}{i-1}}.$$
  We make the obvious modifications in the case of   $F^n_\C$.

% and $[e_1, ie_1]=i[e_1, e_1]=0$. 

\subsection{The  Baker-Campbell-Hausdorff  formula}\label{BCH}

Let  $G$ be a simply connected nilpotent Lie group with Lie algebra $\mathcal{G}$.
  The exponential map  $\text{exp}:  {\mathcal{G}}\ra G$ is a diffeomorphism.  
   One can then
     pull back the group operation from $G$ to get a group stucture 
 on $\mathcal{G}$.     This group structure can be described by the   Baker-Campbell-Hausdorff formula
 (BCH formula in short),  which expresses the product $X*Y$ ($X, Y\in {\mathcal{G}}$)
    in terms of the iterated Lie brackets  of $X$ and $Y$. 
     The group operation in $G$ will be denoted by $\cdot$.  
   The pull-back group operation   $*$  on $\mathcal{G}$ is defined as follows. 
      For $X,   Y\in \mathcal{G}$,   define
   $$X*Y=\text{exp}^{-1}(\text{exp} X\cdot \text{exp} Y).$$
  Then the BCH formula (\cite{CG},    page 11)  says
  $$X*Y=\sum_{n>0}\frac{(-1)^{n+1}}{n}\sum_{p_i+q_i>0, 1\le i\le n}\frac{(\sum^n_{i=1}(p_i+q_i))^{-1}}
  {p_1!q_1!\cdots p_n!q_n!} (\text{ad} X)^{p_1}(\text{ad} Y)^{q_1}\cdots   (\text{ad} X)^{p_n}(\text{ad} Y)^{q_n-1}Y$$
  where $\text{ad}\, A (B)=[A,   B]$.     If  $q_n=0$,  the term in the sum is  $\cdots (\text{ad} X)^{p_n-1}X$;
  if  $q_n>1$  or if $q_n=0$  and $p_n>1$, then the term  is zero. 
  The first a few terms are well known,
\b{align*}
 X*Y& =X+Y+\frac{1}{2}[X,Y]+\frac{1}{12}[X,[X,Y]]
-\frac{1}{12}[Y, [X, Y]]\\
& -\frac{1}{48}[Y,[X,[X,Y]]]-\frac{1}{48}[X,[Y,[X,Y]]]+(\text{commutators  in five or more terms}).
\end{align*}

\b{Le}\label{cbhf}
There are universal constants $c_j\in \Q$,  $2\le  j\le n-1$ with the following property:  
     for any $X\in \mathcal{F}_K^n$ and any $Y\in Ke_2\oplus V_2\oplus \cdots \oplus V_n$, we have
\b{equation}\label{e2.3.10}
X*Y=X+Y+\frac{1}{2}[X, Y]+\sum_{j=2}^{n-1} c_j(\text{ad} \,X)^j \, Y
\end{equation}
  and 
 % (2)   for any $X\in \mathcal{F}_K^n$ and any $Y\in V_2\oplus \cdots \oplus V_n$, we have
\b{equation}\label{e2.3.11}
Y*X=Y+X+\frac{1}{2}[Y, X]+\sum_{j=2}^{n-1} (-1)^j c_j(\text{ad} \,X)^j \, Y.
\end{equation}

\end{Le}

\b{proof}
    (\ref{e2.3.10})     follows from the fact that $[Z_1, Z_2]=0$ for any two 
$Z_1, Z_2\in  Ke_2\oplus   V_2\oplus \cdots \oplus V_n$.   In particular,
  $[Y, [X, Y]]=0$, so the only possible  nonzero terms in the BCH  formula 
  are multiples of $(\text{ad} \,X)^j \, Y.$  We stop at $j=n-1$ since 
   % $Y\in V_2\oplus \cdots \oplus V_n$  and  
  $ \mathcal{F}_K^n$   has  $n$ layers. 

 (\ref{e2.3.11})   follows from (\ref{e2.3.10})     by taking inverse of both sides and then replacing $X$ with $-X$ and $Y$ with $-Y$.

\end{proof}

\b{Cor}\label{cbh}  In $\mathcal{F}_K^n$
  the following holds for all $t,  x_2,  \cdots,   x_{n+1}\in K$:
    $$(-te_1)*\sum_{j=2}^{n+1} x_j e_j *te_1=\sum_{j=2}^{n+1} x'_j e_j,$$
   where $x'_2=x_2$,  
  $x'_j=x_j-tx_{j-1}+ G_j$ for $3\le j\le n+1$.  Here   $G_j$ is a polynomial  of $t$ and the $x_i$'s  
   and each of  its terms 
  has   a factor $t^k$ for some $k\ge 2$.

\end{Cor}

\b{proof}
We apply Lemma \ref{cbhf}   to $Z:=(-te_1)*\sum_{j=2}^{n+1} x_j e_j $  and obtain:
  $$Z=(-te_1)*\sum_{j=2}^{n+1} x_j e_j =-te_1+\sum_{j=2}^{n+1}x_j e_j +\frac{1}{2}\big[-te_1,  \sum_{j=2}^{n+1}x_j e_j\big]+\sum_{j=4}^{n+1}H_j  e_j,$$
  where $H_j$ 
is a polynomial  of $t$ and the $x_i$'s  
   and each of  its terms 
  has   a factor $t^k$ for some $k\ge 2$.
  Notice that  % we only need to use the first 3 terms since
 the last term  in the formula in Lemma \ref{cbhf}   
contains only higher degree terms   in $t$.   Now we apply the BCH formula to $Z* te_1$. The first 3 terms are
 \b{align*}
Z+te_1+\frac{1}{2}[Z, te_1]& =\sum_{j=2}^{n+1}x_j e_j +\frac{1}{2}\big[-te_1,  \sum_{j=2}^{n+1}x_j e_j\big]+\sum_{j=4}^{n+1}H_j  e_j+\frac{1}{2}\big[\sum_{j=2}^{n+1}x_j e_j, te_1\big]+\sum_{j=4}^{n+1}I_j  e_j\\
&=x_2e_2 +\sum_{j=3}^{n+1}(x_j-t x_{j-1})e_j+\sum_{j=4}^{n+1}(H_j+I_j)e_j,
\end{align*}
where $I_j$ 
is a polynomial  of $t$ and the $x_i$'s  
   and each of  its terms 
  has   a factor $t^k$ for some $k\ge 2$.  For the iterated brackets  in  the BCH formula for $Z* te_1$, we only need to consider terms of the form $(ad\, Z)^j te_1$ ($2\le j\le n-1$) since all the other terms involve $te_1$  at least twice and so have higher degree in $t$.   However, a direct calculation shows that $(ad\, Z)^2 te_1$   also has higher degree in $t$. 
  So the same is true for  all the terms  $(ad\, Z)^j te_1$,    $2\le j\le n-1$.   The   Corollary follows. 

\end{proof}

\subsection{Various maps}\label{maps}

Here we recall the definitions of quasi-similarity and quasisymmetric maps.

Let $K\ge 1$ and $C>0$. A bijection $F:X\ra Y$ between two
     metric spaces is called a $(K, C)$-\e{quasi-similarity}  if
\[
   \frac{C}{K}\, d(x,y)\le d(F(x), F(y))\le C\,K\, d(x,y)
\]
for all $x,y \in X$.

Clearly a map is a quasi-similarity if and only if it is biLipschitz.  
 The point here is that often there is control on $K$  but not on $C$.  In this case, the notion of quasi-similarity provides more information about the distortion.   To see this,  just compare 
$(1, 100)$-quasi-similarity and $1000$-biLipschitz.

Let $\eta: [0,\i)\ra [0,\i)$ be a homeomorphism.
    A homeomorphism
%A homeomorphism between metric spaces
$F:X\to Y$ between two metric spaces is
\e{$\eta$-quasisymmetric} if for all distinct triples $x,y,z\in X$,
we have
\[
   \frac{d(F(x), F(y))}{d(F(x), F(z))}\le \eta\left(\frac{d(x,y)}{d(x,z)}\right).
\]
  If $F: X\ra Y$ is an $\eta$-quasisymmetry, then
  $F^{-1}:    Y\ra X$ is an $\eta_1$-quasisymmetry, where
$\eta_1(t)=(\eta^{-1}(t^{-1}))^{-1}$. See \cite{V}, Theorem 6.3.
 A    homeomorphism between metric spaces
  is quasisymmetric if it is $\eta$-quasisymmetric for some $\eta$. 

We remark that 
     quasisymmetric homeomorphisms between general metric spaces
  are quasiconformal. In the case of Carnot groups  (and more
  generally  Loewner  spaces),  a   homeomorphism  is quasisymmetric if and only
  if it is quasiconformal, see \cite{HK}.

 The main result in \cite{BKR} says that a quasiconformal   map  between two proper, 
locally Ahlfors   $Q$-regular  ($Q>1$)   metric spaces  is absolutely 
 continuous on almost every curve.   This result
   applies to quasiconformal maps    $F:  G\ra G$  on  Carnot groups. 

Pansu \cite{P2} proved that a quasisymmetric map 
  $F: G_1\ra G_2$  between two Carnot groups is absolutely continuous. That is, 
  a  measurable set $A\subset G_1$ has measure $0$ if and only if $F(A)$ has measure $0$.

\subsection{Pansu differentiability  theorem }\label{pansud}

First the definition:

\b{Def}\label{pansu-d}
 Let $G$ and $G'$
  be two Carnot groups endowed with Carnot metrics,  and $U\subset G$, $U'\subset G'$ open subsets.
   A map $F: U\ra U'$ is \e{Pansu  differentiable}
    at $x\in U$  if there exists a graded  homomorphism
     $L: G\ra G'$ such that
     $$\lim_{y\ra x}\frac{d(F(x)^{-1}*F(y),\, L(x^{-1}*y))}{d(x,y)}=0.$$
          In this case, the graded  homomorphism
     $L: G\ra G'$  is called the \e{Pansu  differential} of $F$ at $x$, and
     is denoted by   $dF(x)$.
\end{Def}

We have the following  chain rule for Pansu differentials:

\b{Le} \e{(Lemma  3.7 in \cite{CC})}  \label{chain}
   Suppose $F_1: U_1\ra U_2$ is Pansu differentiable at $p$, and $F_2: U_2\ra U_3$ is 
  Pansu differentiable at $F_1(p)$.  Then $F_2\circ F_1$ is Pansu differentiable at $p$  and 
  $d(F_2\circ F_1)(p)=dF_2(F_1(p))\circ dF_1(p)$. 
\end{Le}

Notice that the Pansu differential of the identity map   $U_1\ra U_1$    is the identity isomorphism.
   Hence if $F: U_1\ra U_2$ is bijective, $F$ is Pansu differentiable at $p\in U_1$  and $F^{-1}$ is Pansu differentiable at $F(p)$,  then 
  $dF^{-1}(F(p))=(dF(p))^{-1}$.

The following  result (except the terminology)   is due to  Pansu
[P2].

\b{theorem}\label{pan}
  Let $G, G'$ be   Carnot groups, and $U\subset G$, $U'\subset G'$  open subsets.  
 Let $F: U\ra U'$ be a quasiconformal  map.
   Then $F$ is a.e. Pansu  differentiable. Furthermore, at a.e.
   $x\in   U$, the Pansu  differential $dF(x): G\ra G'$ is a graded
   isomorphism.

\end{theorem}

 In  Theorem \ref{pan}  and the proofs below,  \lq\lq a. e." is with
 respect to the   Lesbegue  measure on $\mathcal{G}=G$.

To simplify the exposition, we introduce the following terminology:
\b{Def}\label{goodpoint}
Let  $F: U\ra U'$ be a quasiconformal map  between open subsets of Carnot groups.  
     A point $p\in U$ is called a \e{good point} (with respect to $F$)  if:\newline
   (1)  $F$ is Pansu differentiable at $p$  and  $dF(p)$ is a    graded isomorphism;\newline
  (2)  $F^{-1}$ is Pansu differentiable at $F(p)$  and  $dF^{-1}(F(p))$ is a    graded isomorphism.

\end{Def}

 It follows from Pansu differentiability theorem  and Pansu's theorem on absolute continuity of quasiconformal maps that a.e. $p\in U$ is a good point.

Let $G$ be a Carnot group  with Lie algebra $\mathcal{G}=V_1\oplus \cdots \oplus V_r$.
  A \e{horizontal line}
     in $G$ is the image of a map $\gamma: \R\ra G$ of the form
 $\gamma(t)=g* tv$ for some $g\in \mathcal{G}$ and some $v\in V_1\backslash\{0\}$. 
  A  non-degenerate closed connected subset of a horizontal line is called a horizontal line segment.

\b{Le}\label{pansu=id}
  Let $F: G \ra G$ be a continuous map on a Carnot group $G$. Suppose $F$ is Pansu differentiable a.e.  and is absolutely continuous on almost all curves.  If   the Pansu differential is a.e. the identity isomorphism,  then   
  $F$ is a left translation  of $G$.

\end{Le}

\b{proof}
After composing $F$ with a left translation  we may assume $F(o)=o$.  
 For each 
  $v\in V_1\backslash\{0\}$,    consider the family of horizontal lines consisting of the left cosets of $\R v$.
   Then  for a.e. horizontal line  $L=g* \R v$ in  this family, the  Pansu differential   of $F$ is  a.e.  on $L$ 
  the identity isomorphism  and $F$ is absolutely continuous on $L$.
      It follows that  $F(L)=F(g)* \R v$ and  %$F|_L: L=g\cdot \R v\ra F(g)\cdot \R v$ is a translation; that is, there is some $c_L\in \R$ such that  
   $F(g* tv)=F(g)* t v$ for all $t\in \R$.
    Since this is true for a.e. $L$ in the family and $F$ is continuous, the same is true for all horizontal lines in the family.     Since the vector $v\in V_1$ is arbitrary, the same is true for all horizontal lines.

  Now let $g\in   G$ be arbitrary. Then there exist a finite sequence of points
  $o=p_0$,   $p_1$,  $\cdots$,   $p_n=g$ and horizontal line segment $\alpha_i$ from $p_{i-1}$ to $p_i$.  
   Now the first paragraph applied to $\alpha_1$ implies that $F$ fixes all points on $\alpha_1$.   
 Now an induction argument shows that $F$ fixes all points on $\alpha_i$.   In particular, $F(g)=g$.

\end{proof}

\section{Quasiconformal maps on the  real  model Filiform groups}\label{s2}

In this Section we will prove Theorem \ref{thA}.  Given a quasiconformal map  $F:  F^n_\R\ra  F^n_\R$  ($n\ge 3$),    we shall pre-compose  and post-compose   $F$ with left translations and graded automorphisms to obtain a map of the form
  $F_h$ described in the Introduction.

\subsection{Graded automorphisms of    $\mathcal{F}^n_\R$}\label{gradedautoreal}

In this   Subsection we will identify all the graded automorphisms of   $\mathcal{F}^n_\R$.

For an element $x\in \mathcal{N}$ in a  Lie algebra, let $\text{rank}(x)$ be the rank of 
 the linear transformation  $ad(x):  \mathcal{N}\ra   \mathcal{N}$,    $ad(x)(y)=[x,y]$.   In other words, 
  $\text{rank}(x)$  is the dimension of the image of $ad(x)$.

\b{Le}\label{l2.1}
Let $\mathcal{F}^n_\R$ be the $n$-step  real  model Filiform algebra.
  Assume $n\ge 3$.   Let $X=ae_1+be_2\in V_1$  be a  nonzero  element in the first layer. 
  Then   $\text{rank}(X)=1$ if and only if $a=0$.

\end{Le}

\b{proof}
Let $Y=y_1e_1+  \cdots   +y_{n+1}e_{n+1}$. Then 
$$[X, Y]=a(y_2e_3+\cdots   + y_ne_{n+1})-by_1e_3.$$ 
    If $a=0$,
  then  $[X,Y]=-by_1e_3$ and so $\text{rank}(X)=1$.  
   If $a\not=0$,   one can vary the $y_i$'s and it is  clear that $\text{rank}(X)=n-1\ge 2$. The   Lemma follows.

\end{proof}

It is clear that a Lie algebra isomorphism preserves the rank of   elements. 
  A graded isomorphism also preserves the first layer  $V_1$.  Hence we obtain:

\b{Le}\label{l2.2}
  Suppose $n\ge 3$.  Then  $h(\R e_2)=\R e_2$ for   every graded isomorphism $h: \mathcal{F}^n_\R
   \ra \mathcal{F}^n_\R$.

\end{Le}

 Given $a_1, a_2\in \R\backslash\{0\}$  and    $b\in \R$,   define
  a    linear map  $h=h_{a_1, a_2, b}:  \mathcal{F}_\R^n\ra  \mathcal{F}_\R^n$  by:
   $$h(e_1)=a_1e_1+be_2, $$
  $$h(e_j)=a_1^{j-2}a_2e_j\;\; \text{for}\;\;  2\le j\le n+1.$$
  It is easy to check that $h$ is a graded isomorphism of $\mathcal{F}_\R^n$.
   The following lemma says that these are the only graded isomorphisms of   $\mathcal{F}_\R^n$.

\b{Le}\label{graded-real}
  A  linear map  $h:  \mathcal{F}_\R^n\ra  \mathcal{F}_\R^n$  is  a graded isomorphism if and only  if 
  $h=h_{a_1, a_2, b}$  for some 
$a_1, a_2\in \R\backslash\{0\}$  and    $b\in \R$.

\end{Le}

\b{proof}
   Let  $h: \mathcal{F}^n_\R
   \ra \mathcal{F}^n_\R$  be a graded isomorphism.     Then  $h({V_i})=V_i$.   
Hence there are constants 
  $a_1,   a_2, b,c \in \R$ and $a_j\in \R$, $3\le  j\le n+1$ such that 
  $h(e_1)=a_1 e_1+ b e_2$,     $h(e_2)=c e_1 +a_2 e_2$,
   $h(e_{j})=a_j e_{j}$ for    $3\le j\le n+1$.   By Lemma \ref{l2.2}  we have 
   $c=0$. 
    Since $h$ is a Lie algebra isomorphism,   for $2\le j\le n$  we have 
$$a_{j+1}e_{j+1}=h(e_{j+1})=h([e_1, e_j])=[h(e_1), h(e_j)]=[a_1e_1+be_2,  a_je_j]=a_1a_je_{j+1},$$
  so $a_{j+1}=a_1a_j$.  It follows that $a_{j}=a_1^{j-2}a_2$ for $2\le j\le n+1$  and $h=h_{a_1, a_2, b}$.

\end{proof}

\b{Le}\label{pdfg}
Let $F: F^n_\R\ra F^n_\R$  be a map. Suppose there are functions
 $f: \R\ra \R$,  $f_j: \R^{n+1} \ra \R$, $2\le j\le n+1$,
  such that 
$$F(x_1e_1*\sum_{j=2}^{n+1} x_j e_j)=f(x_1)e_1*\sum_{j=2}^{n+1} f_j(x) e_j,$$
  where $x=(x_1, \cdots, x_{n+1})$.  
     If  $F$ is Pansu differentiable at $p:=x_1e_1*\sum_{j=2}^{n+1} x_j e_j$  with 
 Pansu differential 
 $dF(p)=h_{a_1(p), a_2(p), b(p)}$,  then 
   $f'(x_1)$ and $\frac{\partial f_2}{\partial x_2}(x)$ exist,  and  
   $$f'(x_1)=a_1(p)\;\;\text{and}\;\; \frac{\partial f_2}{\partial x_2}(x)=a_2(p).$$\newline
      %Furthermore,  
  %   $\frac{\partial f_2}{\partial x_2}(x)$ exists and
% $ \frac{\partial f_2}{\partial x_2}(x)=a_2(x)$.   

\end{Le}

\b{proof}   For $q\in F^n_\R$,  let $y=(y_1, \cdots, y_{n+1})  \in \R^{n+1}$ be determined by  
  $q=y_1e_1*\sum_{j=2}^{n+1} y_j e_j$.     
  By the definition  of Pansu differential,     we have
 \b{equation}\label{e11}
\lim_{q\ra p}\frac{d(F(p)^{-1}* F(q),\, dF(p)(p^{-1}*q))}{d(p,q)}=0.
\end{equation}
   In particular, the above limit is $0$ if $q=p*(y_1-x_1)e_1$  and $y_1\ra x_1$.  
  Let $q=p*(y_1-x_1)e_1$.  Then   $p^{-1}*q=(y_1-x_1)e_1$  and so 
$d(p,q)=d(o,   (y_1-x_1)e_1)=|y_1-x_1|$
  and 
  $$dF(p)( p^{-1}*q)=(y_1-x_1)(a_1(p) e_1+b(p)e_2).$$
Notice that the coefficient of $e_1$ in  $q$ is $y_1$,  and so the coefficient of $e_1$ in
$F(p)^{-1}*F(q)$ is $f(y_1)-f(x_1)$.   The coefficient of $e_1$ in 
 $[-dF(p)(p^{-1}*q)]* F(p)^{-1}* F(q)$ is $f(y_1)-f(x_1)-a_1(p) (y_1-x_1)$.
Hence 
$$d(F(p)^{-1}* F(q), \,dF(p)(p^{-1}*q))\ge  |f(y_1)-f(x_1)-a_1(p) (y_1-x_1)|.$$
 Now (\ref{e11}) implies:
  $$\lim_{y_1\ra x_1}\frac{  |f(y_1)-f(x_1)-a_1(p) (y_1-x_1)|}{|y_1-x_1|}=0.$$
    Hence   $f$ is differentiable  at $x_1$ and $f'(x_1)=a_1(p)$.

  % Suppose   $F$ sends  left cosets   of $\R e_2$ to  left cosets of 
     %$\R e_2$,     
 The proof of the statement about $\frac{\partial f_2}{\partial x_2}(x)$  is similar.  
Let  $q=p*(y_2-x_2) e_2$.    Since $[e_i, e_j]=0$ for $i, j\ge 2$, we can write
  $$q=x_1e_1*\sum_{i=2}^{n+1}x_ie_i*(y_2-x_2) e_2=x_1e_1*(y_2e_2+\sum_{i=3}^{n+1} x_ie_i).$$
  So  $x$ and $y$    differ   only   in the second  coordinate.   
  Notice we have  $p^{-1}* q=(y_2-x_2) e_2$,
   $d(p, q)=|y_2-x_2|$  and  $dF(p)(p^{-1}*q)=a_2(p) (y_2-x_2)e_2$. 
  Also the coefficient of $e_2$ in $F(p)^{-1}*F(q)$    
  is   
  $$f_2(y)-f_2(x)=f_2(x_1, y_2, x_3, \cdots, x_{n+1})-f_2(x_1, x_2, x_3, \cdots, x_{n+1}).$$
Hence 
$$d(F(p)^{-1}* F(q),\,  dF(p)(p^{-1}*q))\ge  |f_2(y)-f_2(x)-a_2(p) (y_2-x_2)|.$$
 Now (\ref{e11}) implies:
  $$\lim_{y_2\ra x_2}\frac{  |f_2(y)-f_2(x)-a_2(p) (y_2-x_2)|}{|y_2-x_2|}=0.$$
   Hence    $\frac{\partial f_2}{\partial x_2}(x)$ exists and
 $ \frac{\partial f_2}{\partial x_2}(x)=a_2(p)$.

\end{proof}

\subsection{Quasiconformal implies biLipschitz}\label{qcbilip1}

In this   Subsection we show that every quasiconformal map of $F^n_\R$  ($n\ge 3$) is biLipschitz.

  Let   $n\ge 3$  and  $F: {F}^n_\R  \ra  {F}^n_\R$   be an  $\eta$-quasisymmetric   map
  for some $\eta$.    
 Pansu's differentiability theorem says that $F$ is Pansu differentiable a.e. and  the Pansu differenntial   is   a.e.   a graded isomorphism.  
Notice that the  Lie subalgebra  generated by $\R e_2$ is itself. 
  We shall abuse notation and also  denote the  corresponding   connected subgroup of 
  $ {F}^n_\R$     by $\R e_2$.
 It follows from Fubini's theorem that for   a.e. left coset
  $L$  of $\R e_2$,  the map $F$ is Pansu differentiable a.e. on $L$  and the Pansu differential is a.e. a graded isomorphism on $L$.    
By Lemma \ref{l2.2}  
  $dF(x)(\R e_2)=\R e_2$ for a.e. $x\in {F}^n_\R$.  
  %Notice that the subgroup generated by $\R e_2$ is itself.  
Now  the following result 
     implies that $F$ sends 
  left cosets of  $\R e_2$ to the left cosets of $\R  e_2$.

\b{Prop}\label{p3} \e{(Proposition   3.4  in \cite{X3}) }
 Let  $G$  and $G'$ be two Carnot groups,  $W\subset V_1$,  $W'\subset V'_1$ be 
  subspaces. Denote by $\mathcal{G}_W\subset \mathcal{G}$ and 
$\mathcal{G}'_{W'}\subset \mathcal{G}'$   respectively the Lie subalgebras
  generated by $W$ and $W'$.    Let $H\subset G$ and $H'\subset G'$ respectively  be the 
  connected   Lie
  subgroups of $G$ and $G'$    corresponding to 
  $\mathcal{G}_W$  and    $\mathcal{G}'_{W'}$.  
 Let $F: G\ra G'$ be a quasisymmetric  homeomorphism.  
  If     $dF(x)(W)\subset W'$ for a.e. $x\in G$,   
     %a.e.  the Pansu-differential of $F$ sends $W$ into $W'$,  
     then  
   $F$ sends left cosets of $H$ into left  cosets of $H'$.  

  \end{Prop}

To simplify the exposition, we introduce the following terminology (see Definition \ref{goodpoint}
  for the  definition of  \lq\lq good point"):

\b{Def}\label{goodcoset}
Let   $n\ge 3$  and  $F: {F}^n_\R  \ra  {F}^n_\R$   be a   quasisymmetric   map.
       A left coset $L$ of $\R e_2$ is called a \e{good left coset}  (with respect to $F$)   if the following conditions hold:\newline
  (1)   a.e. point $p\in L$ is a good point  with respect to $F$;\newline
  (2)  $F|_L$ is absolutely continuous;\newline
  (3)  $F^{-1}|_{F(L)}$ is absolutely continuous.  

\end{Def}

\b{Le}\label{goodlcoset}
Let   $n\ge 3$  and  $F: {F}^n_\R  \ra  {F}^n_\R$   be a   quasisymmetric   map.
  Then a.e. left  coset of  $\R e_2$  is a good left coset.
\end{Le}

\b{proof}
Since a.e. point   is a good point, Fubini's theorem implies  (1)  holds  for a.e. left coset $L$.
  The main theorem of \cite{BKR} implies  $F|_L$  is absolutely continuous on  a.e. left coset  $L$ of  $\R e_2$.
   Since $F^{-1}$ is also quasisymmetric,  $F^{-1}|_{L'}$ is  
absolutely continuous on  a.e. left coset  $L'$ of  $\R e_2$.
  Since quasisymmetric maps preserve the conformal modulus of curve families and $F$ maps left cosets of $\R e_2$ to left cosets of $\R e_2$, we see that  (2) and (3) hold for a. e.  left coset $L$.

\end{proof}

   For a  good  point  $x\in F_\R^n$,    % where $F$ is Pansu differentiable  and the Pansu 
    %differential  $dF(x)$ is a  graded isomorphism,  
 let $b(x),  a_1(x),  a_2(x)\in \R$   be such that 
     $dF(x)=h_{a_1(x),  a_2(x), b(x)}$,      see    Lemma \ref{graded-real}.
    Then  $dF(x)(e_1)=a_1(x) e_1+ b(x) e_2$,  $dF(x)(e_j)= (a_1(x))^{j-2} a_2(x) e_j$  for $2\le  j\le  n+1$.
  Set $a_{n+1}(x)=(a_1(x))^{n-1} a_2(x)$.  
     Notice 
 $dF(x)(v)=a_{n+1}(x)v$   for any $v\in V_n=\R e_{n+1}$.

\b{Le}\label{second}
Let $L$ be a  good  left coset of $\R e_2$.  % where $F$ is Pansu  differentiable a.e.  and 
   % the Pansu 
 % differential   is a  graded isomorphism a.e.   
  Then there is a constant $0\not=a_L\in \R$  such that $a_{n+1}(x)=a_L$  for a.e. $x\in L$.

\end{Le}

\b{proof} 
   Suppose there  exist  good  points    $p, q\in L$ such that 
  %$F$ is 
  %Pansu differentiable   at $p, q$,   $dF(p)$ and $dF(q)$ are graded isomorphisms and 
$a_{n+1}(p)\not=a_{n+1}(q)$.  
 We shall show that this implies $a_{n+1}(z)\ra +\infty$ as $z\ra \infty$ along $L$.  
 This provides a contradiction since the same claim applied to $F^{-1}$  implies $a_{n+1}(z)\ra 0$ as 
$z\ra \infty$ along $L$;  see the remark after Lemma \ref{chain}  about the relation between the Pansu differential of $F$ and that of   $F^{-1}$.

   We may  
 assume   $|a_{n+1}(p)| \ge |a_{n+1}(q)|$.  
  %After pre-composing and post-composing   with left translations we may assume 
  %  $p=F(p)=o$  is the origin.   
   Notice that  $L=p*\R e_2$  and $F(L)=F(p)*\R e_2$.  
   In particular,  there are $y, y_0\in \R$ such that  $q=p*ye_2$ and $F(q)=F(p)*y_0 e_2$.  
  After    composing with the  inverse of   $dF(p)$      we may assume $dF(p)=Id$ is the identity 
  isomorphism.  Then  $a_{n+1}(p)=1$ and  either $|a_{n+1}(q)|<1$  or 
   $a_{n+1}(q)=-1$.  
  %In particular,  there are $y, y_0\in \R$ such that  $q=ye_2$ and $F(q)=y_0 e_2$.  
  We  shall  consider the image of the left coset $r^n  e_{n+1}  * L$ under the map $F$  as $r\ra 0$.

  Notice that $d(o, r^n e_{n+1})=r$.   
Since $dF(p)=Id$,   the definition of Pansu differential implies that  there   exists some    $\tilde{x}=\sum_i \tilde{x}_i e_i$  with $\tilde{x}_i=\tilde{x}_i(r)$ such that 
  $d(o, \tilde{x})=o(r)$  and  $F(r^n e_{n+1}*p)=F(p)*r^n e_{n+1} * \tilde{x}$.   
  Similarly, there is some  $\tilde{y}=\sum_i \tilde{y}_i e_i$  with $\tilde{y}_i=\tilde{y}_i(r, y)$ such that 
  $d(o, \tilde{y})=o(r)$  and  $F(r^n e_{n+1}*q)=F(q)*a_{n+1}(q)r^ne_{n+1}* \tilde{y}$. 
   %y_0e_2*a_3(q) r^2 e_3 * \tilde{y}$.   
  For later use, we notice here that $\tilde{x}_1,  \tilde{x}_2,  \tilde{y}_2=o(r)$  and 
     $\tilde{x}_{n+1},   \tilde{y}_{n+1}=o(r^n)$.    
    Since   $F$ sends left cosets of $\R e_2$ to   left cosets of $\R e_2$,   we have 
 $$L_r:=F(r^n e_{n+1}* L)=F(r^n e_{n+1}*p)* \R e_2=F(p)*r^n e_{n+1} * \tilde{x}*\R e_2.$$  
  In particular,   $F(r^n e_{n+1}*q)= F(p)*r^n e_{n+1} * \tilde{x}* \tilde{s}_1 e_2$ for some $\tilde{s_1}\in \R$.  
  Using   Lemma \ref{cbhf},   
    %BCH formula, 
    we find two expressions for 
$F(r^n e_{n+1}*q)$:
  \b{align*}
 &F(r^n e_{n+1}*q)\\
&=F(q)*a_{n+1}(q)r^ne_{n+1}* \tilde{y}\\
& =
F(p)*y_0 e_2*a_{n+1}(q) r^n e_{n+1} * \tilde{y}\\
& =F(p)* \big(\tilde y_1 e_1+\{\tilde y_2 +y_0\}e_2+S_1+\{\tilde y_{n+1}+a_{n+1}(q) r^n +(-1)^{n-1}c_{n-1}\tilde y_1^{n-1} y_0\}e_{n+1}\big)
\end{align*}
%\tilde{y}_1e_1+ (y_0+\tilde{y}_2)e_2+\big(a_3(q)r^2+\tilde{y}_3-\frac{1}{2}y_0\tilde{y}_1\big)e_3+S_1,
%\end{align*}
  and 
  \b{align*}
F(r^n e_{n+1}*q)& =F(p)*r^n e_{n+1} * \tilde{x}* \tilde{s}_1 e_2\\
  & =F(p)*\left(\tilde{x}_1e_1+\{\tilde{x}_2+\tilde{s}_1\}e_2+S_2+ \{\tilde{x}_{n+1}+r^n+ c_{n-1}\tilde{x}_1^{n-1}\tilde{s}_1\}e_{n+1}\right), 
\end{align*}
  where   $S_1$  and $S_2$ are   linear  combinations    of   the  $e_i$'s  
   with $3\le i\le n$.  
 Comparing the coefficients we obtain:
   $\tilde{y}_1=\tilde{x}_1$,  $\tilde{s}_1=y_0+(\tilde{y}_2-\tilde{x}_2)$ and
 $$a_{n+1}(q)r^n+\tilde{y}_{n+1}+(-1)^{n-1}c_{n-1}\tilde{y}_1^{n-1}y_0=r^n+\tilde{x}_{n+1}+c_{n-1}\tilde{x}_1^{n-1}\tilde{s}_1.$$
    Pluging the first two into the third we obtain:
 \b{equation}\label{e3.10}
(1-a_{n+1}(q))r^n+(1-(-1)^{n-1})c_{n-1}\tilde x_1^{n-1} y_0=(\tilde y_{n+1}-\tilde x_{n+1})-c_{n-1}\tilde x_1^{n-1}(\tilde y_{2}-\tilde x_{2}).
\end{equation}
  %$$y_0\tilde{x}_1=\big(a_3(q)-1\big)r^2+(\tilde{y}_3-\tilde{x}_3)-\frac{1}{2}\tilde{x}_1(\tilde{y}_2-\tilde{x}_2).$$
Since    %$a_{n+1}(q)<1$  and 
$\tilde{x}_1,  \tilde{x}_2,  \tilde{y}_2=o(r)$,    
     $\tilde{x}_{n+1},   \tilde{y}_{n+1}=o(r^n)$,
  we see that we must have $a_{n+1}(q)=1$  if $c_{n-1}=0$  or $n-1$ is even.   
  In these cases  the  Lemma holds. So from now on we shall assume that 
  $n-1$ is odd  and  $c_{n-1}\not=0$. 
Then  (\ref{e3.10}) implies that, 
    for all sufficiently small $r$,   we have 
  %$ y_0\tilde{x}_1<0$,
    %\b{equation}\label{e1}
 %| \tilde{x}_2|\le r/10
%\end{equation} 
  %    and  
\b{equation}\label{e22}
 -\frac{11}{10}(1-a_{n+1}(q))r^n  \le  2 c_{n-1}  \tilde{x}_1^{n-1}  y_0    \le  
   -\frac{9}{10}(1-a_{n+1}(q))r^n .
\end{equation}
  % for all sufficiently small $r$. 

  Our next goal is to bound $|a_{n+1}(z)|$ from below for $z\in L$.  For this 
      we  shall    bound from below  the distance from a point on $F(L)$  to 
$L_r$.  %:=F(r^2 e_3  *\R e_2)$.    %  to  a point on $\R e_2$.    % $F(y)=y_0 e_2$.   
So we fix    $q_2=F(p)*s_2 e_2\in F(L)$  and  let
$q_1\in  L_r$  vary.  % F(r^2 e_3 * \R e_2)$.  %  and $q_2\in \R e_2$. 
     Then 
  $q_1=F(r^ne_{n+1}*p)*s_1 e_2=F(p)*r^n e_{n+1} * \tilde{x}* s_1e_2$     with $s_1\in \R$. 
  Then $d(q_2, q_1)=d(o,  (-s_2e_2)*r^n e_{n+1} * \tilde{x}* s_1e_2)$.  
  Using   Lemma \ref{cbhf} twice we obtain 
    %Using BCH formula one calculates that 
  \b{align*}
  & (-s_2e_2)*r^n e_{n+1} * \tilde{x}* s_1e_2\\
& =\tilde{x}_1 e_1+\{\tilde{x}_2-s_2+s_1\}e_2+S+\{r^n+\tilde{x}_{n+1}+c_{n-1}\tilde{x}_1^{n-1}s_2+
      c_{n-1}\tilde{x}_1^{n-1}s_1\}e_{n+1},
\end{align*}
  where $S$ is a linear   combination    of   the  $e_i$'s    with $3\le i\le n$.  
    Write   $s_2=ty_0$  and $s_1-s_2=u$.   
    Fix  any $t$ with   $|t|\ge \frac{100}{1-a_{n+1}(q)}$.   
   Notice $\frac{100}{1-a_{n+1}(q)}\ge 50$  since $|a_{n+1}(q)|\le 1$.  
 Let $r$ be sufficiently small so that   $| \tilde{x}_2|\le r/10$  
   and 
   (\ref{e22}) holds. 
  If   $s_1$ is such that  $|u|\ge  \sqrt{|t|} r $,    then   %(\ref{e1}) implies  
  $d(q_2, q_1)\ge  |(\tilde{x}_2-s_2+s_1)|  \ge  \sqrt{|t|} r/2 \ge \frac{1}{2}|t|^{\frac{1}{n}}\cdot r$. 
  If   $|u|\le  \sqrt{|t|} r $,    then   (\ref{e22})  and the assumption on 
    $|t|$ imply the following holds for sufficiently small $r$:
\b{align*}
 d(q_2, q_1) & \ge \left|r^n+\tilde{x}_{n+1}+c_{n-1}\tilde{x}_1^{n-1}s_2+
      c_{n-1}\tilde{x}_1^{n-1}s_1\right|^{\frac{1}{n}}\\
& =\left|r^n+\tilde{x}_{n+1}+2c_{n-1}\tilde{x}_1^{n-1}s_2+
c_{n-1}\tilde{x}_1^{n-1}u\right|^{\frac{1}{n}}\\
 & =\left|r^n+\tilde{x}_{n+1}+t\cdot2c_{n-1}\tilde{x}_1^{n-1}y_0+
c_{n-1}\tilde{x}_1^{n-1}u\right|^{\frac{1}{n}}\\
 %=|r^2+\tilde{x}_3+ty_0\tilde{x}_1+\frac{1}{2}u\tilde{x}_1|^{\frac{1}{2}}\\
  &  \ge
\left|\frac{1}{2}\cdot t\cdot2c_{n-1}\tilde{x}_1^{n-1}y_0\right|^{\frac{1}{n}}\\
   %|ty_0\tilde{x}_1/2|^{\frac{1}{2}}\ge |t(1-a_3(q))r^2/4|^{\frac{1}{2}}\\
  &  \ge  \left|\frac{9(1-a_{n+1}(q))}{20} |t| r^n\right|^{\frac{1}{n}}\\
&=\left|\frac{9(1-a_{n+1}(q))}{20}\right|^{\frac{1}{n}} \cdot |t|^{\frac{1}{n}}\cdot r.
%&=\frac{\sqrt{1-a_3(q)}}{2}\sqrt{|t|}r.
\end{align*}
  It follows that
$d(q_2,  L_r)\ge c {|t|^{\frac{1}{n}}}r$, where $c=\min\{\frac{1}{2}, \left|\frac{9(1-a_{n+1}(q))}{20}\right|^{\frac{1}{n}}  \}$.
  This implies that $a_{n+1}(F^{-1}(q_2))\ge c {|t|^{\frac{1}{n}}} \ra +\infty$ as $|t|\ra \infty$, finishing the proof of our claim.

\end{proof}

  Recall that $F$  is $\eta$-quasisymmetric.

\b{Le}\label{new1}
   %There exists  a constant $K_1$ that depends only on $\eta$ with the following property:
      For    every good     left coset  $L$  of $\R e_2$   %where $F$ is Pansu  differentiable a.e.     and   the Pansu 
  %  differential   is a  graded isomorphism  a.e.,   
 the restriction 
  $F|_L:  L  \ra F(L)$  is a $(\eta(1),  {|a_L|^{\frac{1}{n}}})$-quasi-similarity.  

\end{Le}

\b{proof}   % By  the main theorem in \cite{BKR},   for a.e. left coset $L$  of $\R  e_2$,  $F|_L$  is absolutely continuous. 
 Notice that $(L, d)$  is isometric to the real line  and under this identification the usual derivative of the map  $F|_L: L\ra F(L)$  is simply $a_2(p)$.  So we only need to control  $a_2(p)$   for  $p\in L$. 
 % Let  $L$  be a  left coset  $L$  of $\R e_2$  such that 
  %$F|_L$  is absolutely continuous,   
 %$F$ is Pansu  differentiable a.e.   on $L$    and   the Pansu 
  % differential   is a  graded isomorphism  a.e  on $L$.  
%Let $x\in L$  be a point where 
   %$F$ is Pansu  differentiable.  By Lemma \ref{second}, $a_3(x)=a_L$.
   Since $F$ is $\eta$-quasisymmetric,   its Pansu differentials are also $\eta$-quasisymmetric. 
   Notice that $dF(p)(e_2)=a_2(p)e_2$,   $dF(p)(e_{n+1})=a_{n+1}(p) e_{n+1}$.  
   So  $d(o, dF(p)(e_2))=|a_2(p)|$  and $d(o, dF(p)(e_{n+1}))=|a_{n+1}(p)|^{\frac{1}{n}}$.  
  Since   $d(o, e_2)=d(o, e_3)=1$,   
  the quasisymmetric condition now implies 
 $$ \frac{1}{\eta(1)}\cdot  |a_{n+1}(p)|^{\frac{1}{n}}\le   |a_2(p)|\le  {\eta(1)} \cdot  |a_{n+1}(p)|^{\frac{1}{n}}.$$
 By Lemma  \ref{second}  there is a constant $a_L$ such that     $a_{n+1}(p)=a_L$  for a.e.  $p\in L$.
   As  $F|_L: L\ra F(L)$   is a homeomorphism between two lines,   
 either $a_2(p)\ge 0$ for a.e. $p\in L$  or $a_2(p)\le 0$ for a.e. $p\in L$.
 The Lemma now follows  from the fact that $F|_L$ is absolutely continuous.

 % Recall that $a_3(x)=a(x) d(x)$ is the   horizontal  Jacobian of $F$ at $x$.  
  %Since  $F$ is $\eta$-quasisymmetric,   there is a  constant $c$ depending only on $\eta$ such that 
   %$$\frac{|dF(x)(v)|^2}{c^2}  \le a_2(x)=a_L\le  c^2\cdot   |dF(x)(v)|^2$$ 
  %for any unit vector  $v\in V_1$.     Letting $v=e_2$ and using $dF(x)(e_2)=d(x)e_2$ we obtain:
%$\frac{\sqrt{a_L}}{c}\le  d_2(x)\le  c\cdot  \sqrt{a_L}.$
  %Since  this holds for a.e. $x\in L$  and  $F|_L$ is absolutely continuous, the Lemma follows. 

\end{proof}

\b{Le}\label{compare}
  For    every  two  good    left cosets 
    $L,  L'$  of $\R e_2$,  we have 
   $\frac{|a_L|}{2^n\eta(1)^{2n}}\le  |a_{L'}|\le  2^n\eta(1)^{2n} |a_L|$.

\end{Le}

\b{proof}  Let $L, L'$ be two  good   left cosets   of  $\R  e_2$.  %  as in Lemma \ref{new1}.  
      We assume  $|a_{L}|\ge  2^n\eta(1)^{2n}  |a_{L'}|$  and will get a contradiction.
 Fix $g_1 \in L$ and $g_2\in L'$.    By Lemma \ref{new1}   
$$d(F(g_1* te_2), F(g_1))\ge |t| {|a_L|^{\frac{1}{n}}}/\eta(1)\ge  
 2|t| \eta(1) {|a_{L'}|^{\frac{1}{n}}}$$  
   and $d(F(g_2* te_2), F(g_2))\le  |t| \eta(1){|a_{L'}|^{\frac{1}{n}}}$.  
  It follows that 
$$\liminf_{t\ra \infty}\frac{d(F(g_1* te_2),  F(g_2*te_2))}{d(F(g_2),   F(g_2*te_2))}\ge 1.$$
  On the other hand,  it is an easy calculation to show that 
    either    $d(g_1*te_2,   \,g_2*te_2)$  is bounded independent of $t$ if 
       the coefficient of  $e_1$ in $g_1$ and $g_2$ are the same,  otherwise  
  $$d(g_1*te_2, g_2*te_2)\sim \sqrt t\;\;  \text{as}\;\;  t\ra \infty.  $$ 
  %$$d(g_1*te_2, g_2*te_2)\sim \sqrt t\;\;  \text{as}\;\;  t\ra \infty.$$
  Hence 
\b{equation}\label{e3.20} 
\frac{d(g_1*te_2, g_2*te_2)}{ d(g_2, g_2*te_2)}\ra 0\;\; \text{as}\;\; t\ra \infty.
\end{equation}  
    The quasisymmetry condition of $F$ implies 
$$\frac{d(F(g_1* te_2),  F(g_2*te_2))}{d(F(g_2),   F(g_2*te_2))}\ra 0   \;\; \text{as} \;\; t\ra \infty.$$
 So we obtain a contradiction.  Similarly we get a contradiction if  
   $|a_{L'}|\ge  2^n\eta(1)^{2n}  |a_{L}|$.

   %$|a_{L'}|\ge  4\eta(1)^4  |a_{L}|.$  

\end{proof}

%Theorem \ref{thA} follows from the following:  

\b{Le}\label{bilip} Suppose $n\ge 3$.  Then  every $\eta$-quasisymmetric map 
  $F:  F^n_\R  \ra    F^n_\R $ is a $(2\eta(1)^4,   C)$-quasi-similarity for some constant   $C$.

\end{Le}

\b{proof}
  Fix a  good  left coset $L_0$.  %   as in Lemma \ref{new1}.  
   Then  
    $F|_{L_0}$  is a $(\eta(1),  {|a_{L_0}|^{\frac{1}{n}}})$-quasi-similarity.   Lemma \ref{new1} and   Lemma 
 \ref{compare} now   imply   that $F|_L$ is a 
  $(2\eta(1)^3,  {|a_{L_0}|^{\frac{1}{n}}})$-quasi-similarity  for   every good   left coset   $L$   of  $\R e_2$.  
 %Let   $L$ be an arbitrary left coset of 
  %$\R e_2$.  Then   Lemma \ref{new1}  and   Lemma \ref{compare}  imply that $F|_L$ is a 
  %$(2c^3,  \sqrt{a_{L_0}})$-quasi-similarity.  
     Since   %by  Lemma \ref{goodlcoset}
   a.e. left coset of $\R e_2$  is a good left coset (Lemma \ref{goodlcoset})   and 
$F$ is continuous,   $F|_L$ is a 
  $(2\eta(1)^3,  |a_{L_0}|^{\frac{1}{n}})$-quasi-similarity  for   every    left coset   $L$   of  $\R e_2$.  
  %it holds true  for all left cosets of $\R e_2$.   
  Now  
  let $p, q\in F_\R^n$ be two   arbitrary points.
  If $p,q$ lie on the same left coset, then
    $$   \frac{{|a_{L_0}|^{\frac{1}{n}}}}{2\eta(1)^3} \cdot d(p, q)\le   d(F(p), F(q))\le   2\eta(1)^3 {|a_{L_0}|^{\frac{1}{n}}}\cdot d(p, q).$$
    Suppose $p, q$ do not lie on the same left coset.   Pick a point $q'$ such that $p, q'$ lie on the same left coset and $d(p,q')=d(p,q)$.   Then the quasisymmetric condition implies 
$$d(F(p), F(q))\le \eta(1)\cdot d(F(p), F(q'))\le  \eta(1)\cdot 2\eta(1)^3 |a_{L_0}|^{\frac{1}{n}}\cdot d(p,q')=2\eta(1)^4             |a_{L_0}|^{\frac{1}{n}}\cdot d(p,q).$$  
   Now the same argument applied to $F^{-1}$   finishes the proof.

\end{proof}

\b{remark}
The arguments in this  Section can be modified to show a local version 
of  Lemma \ref{bilip}: if $F: U\ra V$ is a quasisymmetric map between two open subsets of  $F^n_\R$ with $n\ge 3$,  then $F$ is locally biLipschitz. That is, every point $x\in U$ 
 has a neighborhood $U'$ such that $F|_{U'}$ is biLipschitz. 
  In the proofs of the above Lemmas it is not necessary to let $|t|\ra \infty$. 
 One can get around this by choosing a sufficiently small neighborhood of $x$. For instance,     for any fixed $t>0$,    the quotient in  (\ref{e3.20})  becomes very small if $g_2$ is sufficiently close to $g_1$.

\end{remark}

\subsection{Quasiconformal maps have special forms}\label{specialform}

In this   Subsection we  show that quasiconformal maps  on   $ F^n_\R $   ($n\ge 3$)   have very special forms. 

%finish the proof of Theorem \ref{thA}.  

Let $n\ge 3$  and  
  $F:  F^n_\R  \ra    F^n_\R $   be  an $\eta$-quasisymmetric map  for some homeomorphism $\eta: [0, \infty)\ra [0, \infty)$.
     By  Lemma \ref{bilip},
        $F$ is $M$-biLipschitz for some $M\ge 1$. 
Let $G=\R{e_2}\oplus V_2\oplus  \cdots \oplus V_n\subset F_\R^n$.  Notice that $G$ is a subgroup of $F_\R^n$.

\b{Le}\label{anotherfo}
  $F$ sends left cosets of $G$ to left cosets of $G$.   
\end{Le}

\b{proof}
 An easy calculation using the BCH formula  shows  that two left cosets of $\R e_2$ lie in the same left coset of $G$ if and only if the Hausdorff distance between them is finite.   Now the Lemma follows   since by    Lemma \ref{bilip}
    $F$ is biLipschitz.  % we have the following:

\end{proof}

\b{Le}\label{fisbilip}
    %Lemma \ref{anotherfo}  implies  
   There is   a    $M$-biLipschitz  homeomorphism $f: \R\ra \R$ such that $F(x_1e_1*G)=f(x_1)e_1*G$.  
\end{Le}

\b{proof}
Lemma \ref{anotherfo}  implies  that
   there is   a    homeomorphism $f: \R\ra \R$ such that $F(x_1e_1*G)=f(x_1)e_1*G$.  
We need to prove  that $f$ is biLipschitz. 
  % Recall that $F(x_1e_1*G)=f(x_1)e_1*G$ for  any   $x_1\in \R$. Also 
 % By  Lemma \ref{bilip},
   %     $F$ is $M$-biLipschitz for some $M\ge 1$. 
  %By   Lemma \ref{anotherfo} $F$ sends left cosets of $G$ to left cosets of $G$. 
      %We claim that $f$ is biLipschitz. Indeed,   
    Let  $x_1, x'_1\in \R$.   Notice that $d(x_1e_1*G,\; x'_1e_1*G)=|x_1-x'_1|$.   Pick  $p\in x_1e_1*G$ and $q\in x'_1e_1*G$ with   $d(p,q)=|x_1-x'_1|$. 
   %d(x_1e_1*G, \;x'_1e_1*G)$.
  Then  
 $$|f(x_1)-f(x'_1)|=d(f(x_1)e_1*G, \;  f(x'_1)e_1*G)\le  d(F(p), F(q))\le M\cdot d(p, q)=M\cdot |x_1-x'_1|.$$
  %The arguments in \cite{SX}
    %imply   that $f$  is a biLipschitz function. 
   A similar argument shows that $f^{-1}$  is    also    $M$-Lipschitz.  

\end{proof}

\b{Le}\label{e2}  For   every good   left coset $L$  of  $\R e_2$,    there exists a constant $0\not=a_{2, L}\in \R$ such that 
 $a_2(p)=a_{2,L}$  for a.e. $p\in L$.     Furthermore, for any $g\in L$,    we have 
$F(g* te_2)=F(g)* a_{2,L} te_2$    for all  $t\in \R$.

 % Let $L$ be a left coset of $\R e_2$ on which $F$ is Pansu differentiable a.e. 
  %Then there exists a constant $d_L$ such that  $d(x)=d_L$ for a.e. $x\in L$.  

\end{Le}

\b{proof}
      For $p\in F^n_\R$,  let $x=(x_1, \cdots, x_{n+1})\in \R^{n+1}$  be determined by
  $p=x_1e_1*\sum_{i=2}^{n+1} x_i e_i$.    
  By   Lemma \ref{second}     % for a.e. left coset $L$  of $\R e_2$,   
   we have $a_1(p) a_2(p)=a_3(p)=a_L$   for a.e. $p\in L$.  
      On the other hand,  
 there is a homeomorphism $f: \R\ra \R$ such that $F(x_1e_1*G)=f(x_1)e_1*G$.
    Lemma  \ref{pdfg}   
 implies that  $a_1(p)=f'(x_1)$.    As all the points on $L$ have the same $e_1$ coefficient,  
  %   It follows that
 $a_1(p)$ is   a.e.  constant along $L$.  Therefore $a_2(p)$ is  a.e. a  constant along $L$.  
  By Lemma \ref{pdfg} again,    $a_2(p)=\frac{\partial f_2}{\partial x_2}(x)$.
    %Notice that $a_2(x)$  is simply the derivative of the map $F|_L: L\ra F(L)$  at $x$.   
   Since $F|_L$ is absolutely continuous,   
%  Since  $F$ is  biLipschitz, its  restriction to  $L$  is also biLipschitz.  
  %Therefore  
  $f_2$ is an affine function of $x_2$ (while the other variables are fixed)  and 
   %$F|_L: L\ra F(L)$   must be an affine map and 
the last statement in the Lemma  holds. 

\end{proof}

%Now Lemma  \ref{e2} implies that 
  %for each left coset $L=g\cdot \R e_2$ as in the lemma,  
   %we have  $F(g* te_2)=F(g)* d_L te_2$. 

\b{Le}\label{same}
  There exists a constant $0\not=a_2\in \R$ such that 
  $F(g* te_2)=F(g)*  a_2te_2$ for every $g\in F^n_\R$.

\end{Le}

\b{proof}   Since $F$ is continuous, it suffices to show $a_{2,L}=a_{2, L'}$  for   any    two   good   left cosets 
   $L,  L'$  of $\R e_2$.
  % where $F$ is Pansu  differentiable a.e.
    Suppose  there are  two good  left  cosets $L, L'$  of $\R e_2$   % as in Lemma \ref{e2}
 such that $a_{2,L}\not=a_{2, L'}$.  
 Fix $g_1 \in L$ and $g_2\in L'$.  Then
  $F(g_1*te_2)=F(g_1)* a_{2,L}te_2$  and $F(g_2*te_2)=F(g_2)*a_{2, L'}te_2$.
  Now an easy calculation shows that   either    $d(g_1*te_2, g_2*te_2)$  is bounded independent of $t$ if 
 $L, L'$ lie in the same left coset of $G$,  or  
    $$d(g_1*te_2, g_2*te_2)\sim \sqrt t\;\;  \text{as}\;\;  t\ra \infty$$ 
  if $L, L'$   lie in distinct   left cosets of $G$.  On the other hand,   the  triangle inequality  implies  
 $$d(F(g_1)* a_{2, L}  te_2, \;   F(g_2)*a_{2, L'}te_2)\sim
   |a_{2,L}-a_{2, L'}|t\;\; \text{ as}\;\;   t\ra \infty. $$     
   This contradicts the assumption that $F$ is a quasisymmetric map since
  $$\frac{d(g_1*te_2, g_2*te_2)}{d(g_1*te_2, o)}\ra 0$$ 
 and 
$$\frac{d(F(g_1*t e_2),\;  F(g_2 * t e_2))}{d( F(g_1* t e_2), F(o))}=
\frac{d(F(g_1)* a_{2,L} te_2, \;   F(g_2)*a_{2, L'}te_2)}{d(F(g_1)* a_{2, L} te_2, \;   F(o))}\ra \frac{|a_{2, L}-a_{2, L'}|}{|a_{2,L}|}\not=0. $$

\end{proof}

       By replacing   $F$ with 
    $h_{1,  a_2^{-1}, 0}\circ F$   
we  may assume $a_2=1$.

 Notice that the left cosets  of $G$ with the  metric  induced   from $d$   are    isometric   
      to 
 $\R^n$ with the metric 
 $D((x_j), (y_j))=\sum_j |x_j-y_j|^{\frac{1}{j}}$.  
  It was proved in   Section 15 of  \cite{T}
    that each quasisymmetric map $(\R^n, D)\ra (\R^n, D)$ preserves the foliation
    consisting of    affine subspaces parallel to $\R^i\times \{0\}$  for each  $1\le i \le n-1$.  
  This implies that    there exist    continuous  functions $f_j:=f_j(x_1, x_j, x_{j+1}, \cdots,  x_{n+1})$,  $2\le j\le n+1$,
  such that    $F$ has   the following form:
  \b{equation}\label{e3}
F\big(x_1e_1*\big(\sum_{j=2}^{n+1}x_j e_j\big)\big)=f(x_1)e_1 *\sum_{j=2}^{n+1} f_j\, e_j.
 \end{equation}

Let $E\subset \R$ be the subset consisting   of    
    all $x_1\in \R$ with the following properties:\newline
 (1) $f$  is differentiable at $x_1$; \newline  %and $f^{-1}$ are 
  (2)  almost every point  $p$ in  the   left coset $x_1 e_1* G$  is a good point with respect to $F$.\newline
      By   Fubini's theorem and the fact  that  $f$ is biLipschitz,   we see    that $E$ has full measure in $\R$.

\b{Le}\label{add1}
       For $x_1\in E$,   
%For a.e. $x_1\in \R$,     
 there exist   functions $h_j:=h_j(x_1, x_{j+1}, \cdots, x_{n+1})$,  $2\le   j\le  n+1$,   such that  
  the following holds when $x_1\in E$:
  % for  a.e.  $x_1\in \R$,   $F$ is given by 
   $$F(x_1e_1*(\sum_{j=2}^{n+1}x_j e_j))=f(x_1)e_1 *\sum_{j=2}^{n+1} \{(f'(x_1))^{j-2}x_j+h_j)\}e_j.$$
\end{Le}

\b{proof}
  %By  Lemma \ref{bilip}  $F$ is $M$-biLipschitz for some $M\ge 1$.  
 We shall first show that 
$f_j=f_j(x_1, x_j, x_{j+1}, \cdots,  x_{n+1})$,  $2\le j\le n+1$,  is  Lipschitz  in  $x_j$. 
     Let 
  $$p=x_1e_1*(x_2e_2+\cdots +x_{n+1} e_{n+1})$$  
and 
  $$q=x_1e_1*(x_2 e_2+\cdots +  x_{j-1}e_{j-1}+
 y_je_j+x_{j+1}e_{j+1}+\cdots +x_{n+1}e_{n+1}).  $$
    Notice the only difference between $p$ and $q$   is   in the coefficient of $e_j$.  
  We have $(-p)*q=(y_j-x_j)e_j$  and    
       $d(p, q)=|y_j-x_j|^{\frac{1}{j-1}}$.
  %=d(o, -p*q)=d(o,  (y_j-x_j)e_j)=|y_j-x_j|^{\frac{1}{j}}$.
  %By  Lemma \ref{bilip}  $F$ is $M$-biLipschitz for some $M\ge 1$.  
 Hence 
  $$d(F(p), F(q))\le M\, |y_j-x_j|^{\frac{1}{j-1}}.  $$
On the other hand,   by (\ref{e3})  the coefficient  of $e_j$ in 
    $(-F(p))*F(q)$ is:
 \b{equation}\label{e6}
f_j(x_1, y_j, x_{j+1}, \cdots, x_{n+1})- f_j(x_1, x_j, x_{j+1}, \cdots, x_{n+1}).
\end{equation}
Hence
 $$|f_j(x_1, y_j, x_{j+1}, \cdots, x_{n+1})- f_j(x_1, x_j, x_{j+1}, \cdots, x_{n+1})|^{\frac{1}{j-1}}\le d(F(p), F(q))\le M\cdot  |y_j-x_j|^{\frac{1}{j-1}},  $$
  and so 
$$|f_j(x_1, y_j, x_{j+1}, \cdots, x_{n+1})- f_j(x_1, x_j, x_{j+1}, \cdots, x_{n+1})|\le   M^{j-1} |y_j-x_j|.$$
  Hence  $f_j$ is Lipschitz in $x_j$.  %   in particular, it is absolutely continuous in $x_j$.  

Set  $G_j=\R e_2\oplus V_2\oplus \cdots \oplus V_{j-2}\oplus V_{j}\oplus \cdots \oplus  V_n$.  
   For   $x_1\in E$,  define a subset $E_{x_1}\subset   G_j$ as follows:
 $$w=x_2e_2+  \cdots  + x_{j-1}e_{j-1}+  x_{j+1}e_{j+1}+ \cdots+ x_{n+1}e_{n+1}\in E_{x_1}$$   
%   $w=(x_2, \cdots, x_{j-1}, x_{j+1}, \cdots, x_{n+1})\in E_{x_1}$ 
 if and only if 
  the point  %$F$ has nonsingular  Pansu differential   at 
 $$p=p(x_1, w, x_j):=x_1e_1*(x_2e_2+  \cdots  + x_{j-1}e_{j-1}+ x_je_j+ x_{j+1}e_{j+1}+ \cdots+ x_{n+1}e_{n+1})$$
      is   a  good  point    for a.e. $x_j\in \R$.  
  Since   a.e. point   in  $ x_1 e_1*G$ is a good point,    
%Since $f(x_1)$  is differentiable at a.e.  $x_1\in \R$,  and  $F$  is Pansu differentiable a.e.  
%and its Pansu differential is a.e. a graded isomorphism,  
   Fubini's theorem  implies that        $E_{x_1}$ has full measure in   $G_j$.  
   % for 
 %  $x_1\in \R$,        and a.e.  
 %$$x_2e_2+  \cdots  + x_{j-1}e_{j-1}+  x_{j+1}e_{j+1}+ \cdots+ x_{n+1}e_{n+1}\in   \R e_2\oplus V_2\oplus \cdots \oplus V_{j-2}\oplus V_{j}\oplus \cdots \oplus  V_n,$$   
  %the point  %$F$ has nonsingular  Pansu differential   at 
 %$$p=x_1e_1*(x_2e_2+  \cdots  + x_{j-1}e_{j-1}+ x_je_j+ x_{j+1}e_{j+1}+ \cdots+ x_{n+1}e_{n+1})$$
   %   is   a  good  point    for a.e. $x_j\in \R$.  
  Fix $x_1\in E$ and $w\in E_{x_1}$.    Let  $x_j$  be such that 
  $p=p(x_1,   w, x_j)$   is   a  good point. 
  %    Fix such a point $p$.    
 By Lemma \ref{pdfg},  
%As seen in the proof of Lemma \ref{e2},
  %  we have 
$a_1(p)=f'(x_1)$.  By our normalization $a_2=1$, we obtain 
   $a_j(p)=(f'(x_1))^{j-2}$.   Let $q$ be as defined  at the beginning of the proof.  
    %It follows that 
  Then $dF(p)((-p)*q)=(f'(x_1))^{j-2}(y_j-x_j)e_j$.  
   By the definition of   Pansu differential, we have
 \b{equation}\label{e5}
\frac{d(o, \;   [-dF(p)((-p)*q)]*(-F(p))*F(q))}{d(p, q)}\ra 0
\end{equation}
  as $y_j\ra x_j$.  
  Notice that $(-F(p))*F(q)\in G$ and the coefficient of $e_j$ in $(-F(p))*F(q)$ 
  is  given by (\ref{e6}).    It follows that the coefficient of $e_j$ in 
  $[-dF(p)((-p)*q)]*(-F(p))*F(q)$  is 
 $$A_j:=f_j(x_1, y_j, x_{j+1}, \cdots, x_{n+1})- f_j(x_1, x_j, x_{j+1}, \cdots, x_{n+1})-(f'(x_1))^{j-2}(y_j-x_j).$$
  Now  (\ref{e5}) implies 
$$\frac{|A_j|^{\frac{1}{j-1}}}{|y_j-x_j|^{\frac{1}{j-1}}}=\frac{|A_j|^{\frac{1}{j-1}}}{d(p, q)}\le \frac{d(o,  [-dF(p)((-p)*q)]*(-F(p))*F(q))}{d(p, q)}\ra 0.$$
  It follows that $\frac{\partial f_j}{\partial x_j}(x)=(f'(x_1))^{j-2}$. 
  Since this is true for a.e. $x_j$ and   $f_j$ is Lipschitz in $x_j$,  we see that $f_j$ is an affine function of $x_j$ (when the other variables are fixed).  Hence 
  there is some real  number  $H_j:=H_j(x_1, x_{j+1}, \cdots,  x_{n+1})$ such that $f_j=(f'(x_1))^{j-2}  x_j +H_j$.    So far,    $H_j(x_1, x_{j+1}, \cdots,  x_{n+1})$ is  defined   only 
    for those $(x_{j+1}, \cdots, x_{n+1})$ such that  there are $(x_2, \cdots, x_{j-1})$ with
  $$x_2e_2+\cdots   +x_{j-1}e_{j-1}+x_{j+1}e_{j+1}+\cdots  +x_{n+1}e_{n+1}\in E_{x_1}. $$
  Since $E_{x_1}$ has full measure, Fubini's theorem implies that $H_j$   is defined for a.e. 
  $(x_{j+1}, \cdots, x_{n+1})$.

Let $x_1\in E$ be fixed.  
  Set    $$h_j(x_1,   x_j,  x_{j+1}, \cdots,  x_{n+1})=f_j(x_1, x_j, \cdots, x_{n+1})-(f'(x_1))^{j-2}  x_j.$$
    To complete the proof of the Lemma, it 
    suffices to show that $h_j$ is independent of $x_j$.  
  What we proved in the preceding paragraph is:  for a.e. 
 $(x_{j+1}, \cdots, x_{n+1})$,  the function  $h_j(x_1,   x_j,  x_{j+1}, \cdots,  x_{n+1})$  is independent of $x_j$.  Since $h_j$ is continuous  in   $(x_j, x_{j+1}, \cdots, x_{n+1})$,  it is independent of $x_j$ for all 
 $(x_{j+1}, \cdots, x_{n+1})$.  
   
\end{proof}

  Next 
we shall first show that $f$ is an affine function and  then show that $h_i(x_1, x_{i+1}, \cdots, x_{n+1})$ depends only on $x_1$.

Let $x_1\in   E$.  % be as in Lemma \ref{add1}.
Let $p=x_1e_1* \sum_{i=2}^{n+1}x_ie_i$  and $q=x_1e_1* (\sum_{i=2}^{n}x_ie_i+\tilde x_{n+1}e_{n+1})$.  Then 
  $d(p,q)=|x_{n+1}-\tilde x_{n+1}|^{\frac{1}{n}}$.    Since  $F$ is $M$-biLipschitz, 
   we have 
  $$d(F(p), F(q))\le M |x_{n+1}-\tilde x_{n+1}|^{\frac{1}{n}}.$$  
  On the other hand,    
\b{align*}
& ( -F(q))*F(p)\\
 & =\sum_{i=2}^n(h_i(x_1, x_{i+1}, \cdots,  x_{n+1})-h_i(x_1, x_{i+1}, \cdots,  \tilde x_{n+1}))e_i +f'(x_1)^{n-1}(x_{n+1}-\tilde x_{n+1})e_{n+1}.  
\end{align*}
  It follows that  for $2\le i\le n$,  
\b{equation}\label{e8}
|h_i(x_1, x_{i+1}, \cdots,  x_{n+1})-h_i(x_1, x_{i+1}, \cdots, \tilde x_{n+1})|\le M^{i-1}\cdot |x_{n+1}-\tilde x_{n+1}|^{\frac{i-1}{n}}.
\end{equation}

\b{Le}\label{f'}
  The   function   $f$  is affine.

\end{Le}

\b{proof}  
   Let $x_1, \tilde x_1\in   E$.  %\R$ be   as in Lemma  \ref{add1}.  
We shall prove that $f'(x_1)=f'(\tilde x_1)$.   The Lemma then follows since $f$ is  biLipschitz and $E$ has full measure in $\R$. 
  Let $p=x_1e_1*x_{n+1}e_{n+1}$  and $q=\tilde x_1e_1*x_{n+1}e_{n+1}$.
      Eventually we will let $x_{n+1}\ra \infty$.  
    Set  $s_{n+1}=f'(x_1)^{n-1}x_{n+1}+ h_{n+1}(x_1)$  and 
   $\tilde s_{n+1}=f'(\tilde x_1)^{n-1}x_{n+1}+ h_{n+1}(\tilde x_1)$. 
    By  Lemma   \ref{add1}
  $$F(p)=f(x_1)e_1*\{(\sum_{i=2}^n h_ie_i)+ s_{n+1}e_{n+1}\}$$
  and 
  $$F(q)=f(\tilde x_1)e_1*\{(\sum_{i=2}^n \tilde h_ie_i)+  \tilde s_{n+1}e_{n+1}\}, $$
  where $h_i=h_i(x_1, 0, \cdots, 0,  x_{n+1})$  and $\tilde h_i=h_i(\tilde x_1, 0, \cdots,  0,   x_{n+1})$.  
  Denote $a:=f(x_1)-f(\tilde x_1)$.  
  Now 
$$(-F(q))*F(p)=\{-(\sum_{i=2}^n  \tilde h_ie_i)-  \tilde s_{n+1}e_{n+1}\}*  a e_1*\{(\sum_{i=2}^n h_ie_i)+ s_{n+1}e_{n+1}\}$$
    Notice that $d(p,q)=|x_1-\tilde x_1|$.  Hence $d(F(p), F(q))\le M\cdot |x_1-\tilde x_1|$  is bounded from above  by a constant independnt of $x_{n+1}$.  
  By    using  Lemma \ref{cbhf}  twice  
   we  find that the coefficient of $e_{n+1}$ in 
   $(-F(q))*F(p)$  is given by:
    \b{align*}
&  s_{n+1}-\tilde s_{n+1}+\frac{1}{2}a (h_n+\tilde h_n)+\sum_{j=2}^{n-1}c_j a^jh_{n-j+1}-
\sum_{j=2}^{n-1}d_j a^j\tilde h_{n-j+1}\\
& =[f'(x_1)^{n-1}-f'(\tilde x_1)^{n-1}]x_{n+1}+h_{n+1}-\tilde h_{n+1}++\frac{1}{2}a (h_n+\tilde h_n)+\sum_{j=2}^{n-1}c_j a^jh_{n-j+1}-
\sum_{j=2}^{n-1}d_j a^j\tilde h_{n-j+1}
%\sum_{i=2}^n c_i [f(x_1)-f(x'_1)]^{n+1-i} (h_i+h'_i) ,  
  \end{align*}
   which is bounded from above by a constant independent of $x_{n+1}$.    
      By  (\ref{e8})
  $$|h_i(x_1, 0, \cdots, 0, x_{n+1})-h_i(x_1, 0, \cdots, 0)|\le M^{i-1}\cdot |x_{n+1}|^{\frac{i-1}{n}},$$
  we see that $|h_i|$ is  bounded from above by a subliner function of $x_{n+1}$  as $x_{n+1}\ra \infty$.  The same is true for $\tilde h_i$.   It follows that  
$[f'(x_1)^{n-1}-f'(\tilde x_1)^{n-1}]x_{n+1}$ is also bounded above by a sublinear function  of $x_{n+1}$. This can happen only when $(f'(x_1))^{n-1}=(f'(\tilde x_1))^{n-1}$.  
Since $f:    \R\ra \R$ is a homeomorphism,   $f'(x_1)$  and $f'(\tilde x_1)$ have the same sign and hence
  $f'(x_1)=f'(\tilde x_1)$.

\end{proof}

Hence there are constants $a\not=0, b\in \R$ such that $f(x_1)=a x_1+b$.  
  After replacing $F$ with $h_{a^{-1}, 1,0}\circ L_{-F(o)}\circ F$,   
  %After composing with a left translation and a graded isomorphism,
 we may assume
  that   $f(x_1)=x_1$ is the identity map.  %  and $F(o)=o$.  
 So  for $x_1\in E$   we have  
  \b{equation}\label{e30}
f_i(x_1, x_{i}, \cdots, x_{n+1})=x_i+ h_i(x_1,  x_{i+1}, \cdots, x_{n+1}).
\end{equation}
  Now we can extend the definition of  $h_i$ to the case when $x_1\notin E$.
    For any  $(x_1, x_i, \cdots, x_{n+1})$,   set 
 $$H_i(x_1, x_i, x_{i+1}, \cdots, x_{n+1})=f_i(x_1, x_i, \cdots, x_{n+1})-x_i.$$ 
  Notice that $H_i$ is continuous in all variables since $f_i$ is. 
    The equality  (\ref{e30}) implies that $H_i(x_1, x_i,\cdots, x_{n+1})$ is independent of $x_i$  when $x_1\in E$.    Since $E$ has full measure  and $H_i$ is continuous,   we conclude that $H_i$ is independent of $x_i$ for all  $x_1\in \R$.    Hence $H_i(x_1, x_i, \cdots, x_{n+1})$ is a function of $x_1$, $x_{i+1}, \cdots, x_{n+1}$  only.      So   we can define 
$$h_i(x_1, x_{i+1}, \cdots, x_{n+1})=H_i(x_1, x_i, \cdots, x_{n+1})$$
   for any  $(x_1,  x_{i+1}, \cdots, x_{n+1})$.  
     Now     the following holds for all points in $F^n_\R$:
  %the map $F$ takes the form:
  $$F(x_1e_1*(\sum_{i=2}^{n+1}x_i e_i))=x_1e_1 *\sum_{i=2}^{n+1} \{(x_i+h_i(x_1,   x_{i+1}, \cdots, x_{n+1})\}e_i.$$

Next we shall  show that  $h_i(x_1, x_{i+1}, \cdots, x_{n+1})$ depends only on $x_1$.  

\b{Le}\label{hi}
  For each $2\le i\le n$, the function $h_i(x_1, x_{i+1}, \cdots, x_{n+1})$ depends only on $x_1$.

\end{Le}

\b{proof}
The idea is very simple: $F$ sends horizontal vectors to horizontal  vectors.
   Notice that for any $g\in F^n_\R$,   the tangent vectors of $g*te_1$  ($t\in \R$)  are horizontal. 
Since $F$ is Pansu differentiable a.e.,   we see that for a.e. $g\in F^n_\R$,  the tangent vector
   of  the curve $F(g*te_1)$ at $t=0$ exists and is horizontal.     We shall calculate this tangent vector.

 Let $g=x_1e_1*\sum_{i=2}^{n+1}x_ie_i$   be a point where $F$ is Pansu differentiable.  
  By   Corollary    \ref{cbh}   we   have:
 $$g*te_1=x_1e_1*\sum_{i=2}^{n+1}x_ie_i*te_1=(x_1+t)e_1*(-te_1)*\sum_{i=2}^{n+1}x_ie_i*te_1
=(x_1+t)e_1*\sum_{i=2}^{n+1} x'_i e_i,$$
   where $x'_2=x_2$  and    %$x'_3=x_3+c_3tx_2$  and 
  $x'_i=x_i-tx_{i-1}+ G_i$ for $3\le i\le n+1$.  Here   $G_i$ is a polynomial of $t$ and the $x_j$'s   and 
  each of  its terms
  has  a  factor  $t^k$  for some $k\ge 2$.      Denote $x'_1=x_1+t$,       
$h_i=h_i(x_1, x_{i+1}, \cdots,  x_{n+1})$  and 
 $\tilde h_i=h_i(x'_1,   x'_{i+1}, \cdots, x'_{n+1})$.   
  Now 
    $$F(g*te_1)=
F(x'_1e_1*(\sum_{i=2}^{n+1}x'_i e_i))=x'_1e_1 *\sum_{i=2}^{n+1} (x'_i+\tilde h_i) e_i.$$
    Since 
$$F(g)=F(x_1e_1*(\sum_{i=2}^{n+1}x_i e_i))=x_1e_1 *\sum_{i=2}^{n+1} (x_i+h_i)e_i,$$
   by using     Corollary  \ref{cbh}   we  obtain  
 \b{align*}
 & (-F(g))*F(g*te_1)   \\ &=\sum_{i=2}^{n+1} (-x_i-h_i)e_i*(te_1) *\sum_{i=2}^{n+1} (x'_i+\tilde h_i) e_i   \\ &
=(te_1)*(-te_1)*\sum_{i=2}^{n+1} (-x_i-h_i)e_i*(te_1) *\sum_{i=2}^{n+1} (x'_i+\tilde h_i) e_i\\
&  =(te_1)*\sum_{i=2}^{n+1}  x''_ie_i* \sum_{i=2}^{n+1} (x'_i+\tilde  h_i)  e_i\\
&=(te_1)*\sum_{i=2}^{n+1} (x''_i+x'_i+\tilde  h_i)e_i,
\end{align*}
  where $x''_2=-[x_2+h_2]$,    
  %  $x''_3=-[x_3+h_3]-c_3t[x_2+h_2]$,  
$x''_j=-[x_j+h_j]+t[x_{j-1}+h_{j-1}]+H_j$ for $3\le j\le n+1$.
   %  and $h_j=h_j(x_1, x_{j+1}, \cdots,  x_{n+1})$.   
 Here   $H_j$ is a polynomial  and all  its terms
  have   degree at least 2 in $t$.   Set    %$h'_i=h_i(x'_1,   x'_{i+1}, \cdots, x'_{n+1})$  and 
      $\tilde{x}_i=x''_i+x'_i+\tilde  h_i$  for $2\le i\le n+1$.  
    Observe that 
  $\tilde{x}_2=-[x_2+h_2]+x_2+\tilde h_2=\tilde h_2-h_2$  and   for $3\le i\le n+1$, 
  $$ \tilde{x}_i=x''_i+x'_i+\tilde h_i=-[x_i+h_i]+t[x_{i-1}+h_{i-1}]+H_i+x_i-tx_{i-1}+G_i+\tilde h_i=\tilde h_i-h_i+th_{i-1}+G_i+H_i.  $$
      By Lemma \ref{cbhf}
 \b{align*}
 & (-F(g))*F(g*te_1)   \\ &=(te_1)*\sum_{i=2}^{n+1}\tilde{x}_ie_i=te_1+\tilde{x}_2e_2+\sum_{i=3}^{n+1}(\tilde{x}_i+\frac{1}{2} t\tilde{x}_{i-1}+I_i)e_i,
\end{align*}
  where $I_i$ is a  polynomial in $t$ and the $\tilde{x}_j$'s and has  degree  at least 2 in $t$.  
  So  % the coefficient of $e_3$ is 
 % $$\tilde{x}_3+\frac{1}{2}t\tilde{x}_{2}+I_3=h'_3-h_3-c_3th_2+G_3+H_3+\frac{1}{2}t(h'_2-h_2)+I_3=h'_3-h_3-c_3th_2+\frac{1}{2}t(h'_2-h_2)+(G_3+H_3+I_3).$$
    for   $3\le i\le n+1$     the coefficient of $e_i$  is:
  \b{align*}
&  \tilde{x}_i+\frac{1}{2}t\tilde{x}_{i-1}+I_i  \\
& =\tilde h_i-h_i+th_{i-1}+G_i+H_i+\frac{1}{2}t [\tilde h_{i-1}-h_{i-1}+th_{i-2}+G_{i-1}+H_{i-1}]+I_i\\
&=\tilde h_i-h_i+th_{i-1} +\frac{1}{2}t(\tilde h_{i-1}-h_{i-1})+J_i,
  \end{align*}
  where $J_i$  is a  polynomial  and has  degree  at least 2 in $t$.  
    Now the fact that the tangent vector of  the curve $F(g*te_1)$ at $t=0$ is horizontal implies that 
  for all $3\le i\le n+1$,  
 $$\lim_{t\ra 0}\frac{ \tilde{x}_i+\frac{1}{2}t\tilde{x}_{i-1}+I_i }{t}=0.$$
 Hence 
$$\lim_{t\ra 0}\frac{\tilde h_i-h_i+th_{i-1} +\frac{1}{2}t(\tilde h_{i-1}-h_{i-1})+J_i}{t}=0.$$
  Clearly $\lim_{t\ra 0}\frac{J_i}{t}=0$. Since $\lim_{t\ra 0}\tilde h_{i-1}=h_{i-1}$, we have 
$\lim_{t\ra 0}\frac{t(\tilde h_{i-1}-h_{i-1})}{t}=0$.  Hence we have 
  \b{equation}\label{e31}
-h_{i-1}=\lim_{t\ra 0}\frac{\tilde h_i-h_i}{t}=\lim_{t\ra 0}\frac{h_i(x'_1, x'_{i+1}, \cdots,  x'_{n+1})-h_i(x_1, x_{i+1}, \cdots,  x_{n+1})}{t}.
\end{equation}
  For $i=n+1$, we have 
  $$-h_n(x_1, x_{n+1})=\lim_{t\ra 0}\frac{h_{n+1}(x'_1)-h_{n+1}(x_1)}{t}=\lim_{t\ra 0}\frac{h_{n+1}(x_1+t)-h_{n+1}(x_1)}{t}=h'_{n+1}(x_1).$$

  What we have shown is that  $-h_n(x_1, x_{n+1})=h'_{n+1}(x_1)$  at every point $x_1e_1*\sum_{i=2}^{n+1}x_i e_i$   where $F$ is Pansu differentiable.    Since $F$ is Pansu differentiable  a.e., 
   Fubini's theorem implies that  for a.e. $x_1\in \R$,  the equality     
 $-h_n(x_1, x_{n+1})=h'_{n+1}(x_1)$   holds for a.e. $x_{n+1}$.  The continuity of $h_n$ implies that 
   $-h_n(x_1, x_{n+1})=h'_{n+1}(x_1)$   holds for   all $x_{n+1}$.     In particular,  for  a.e. $x_1\in \R$,        $h_n(x_1, x_{n+1}) $ is   independent of $x_{n+1}$. Now the continuity of $h_n$    implies that 
   for all $x_1$,   $h_n(x_1, x_{n+1}) $ is   independent of $x_{n+1}$. 
    This shows that $h_n(x_1, x_{n+1})$ is a function of $x_1$ only. 
  Now  an  induction   argument on $i$   using (\ref{e31})  implies that  for   all $2\le i\le n+1$,   
 $h_i$ is a function of $x_1$ only.

\end{proof}

\subsection{Completing the proof  of Theorem \ref{thA}}\label{provethA}

   Here we  finish the proof of Theorem \ref{thA}.

   We shall   first  show that every map of the form $F_h$ is biLipschitz, hence quasiconformal.
  For this part we shall use the Carnot metric $d_c$.     Recall that $d_c$ and $d$ are   biLipschitz equivalent.  
  Let $h$ be   $M$-Lipschitz  and $F:=F_h$ be defined as in the Introduction.   Clearly   $F$   isometrically maps each  left coset of $G$ to   itself.  
Since $h_j(x)=-\int_0^x h_{j-1}(s)ds$, 
   the calculation in 
the proof of  Lemma \ref{hi}  shows that the curve 
  $F(g*te_1)$ is horizontal.    Furthermore,   %the calculation in that proof shows that 
  the coefficients of $e_1$ and $e_2$ in $-F(g)*F(g*te_1)$ are $t$ and $\tilde{x}_2=h_2(x_1+t)-h_2(t)=h(x_1+t)-h(x_1)$.
So the tangent vector of $F(g*te_1)$ is   $e_1+h'(x_1)e_2$, and has length $\le \sqrt{1+M^2}$.
 It follows that for each horizontal line segment $S$ contained in some left coset of $\R e_1$,
   its image  $F(S)$ has length at most $\sqrt{1+M^2}\cdot \text{length}(S)$.
    Now  let $p, q\in F^n_\R$ be arbitrary.   Then $p\in x_1e_1*G$ and $q\in x'_1e_1*G$ for some 
  $x_1, x'_1\in \R$.   If 
$$d_c(q, p*G)\ge \frac{1}{10\sqrt{1+M^2}} d_c(p,q),$$     then 
  $$d_c(F(q), F(p))\ge d_c(F(q*G), F(p*G))=d_c(q*G, p*G)=d_c(q, p*G)\ge \frac{1}{10\sqrt{1+M^2}} d_c(p,q). $$
    Now suppose   $$d_c(q, p*G)\le \frac{1}{10\sqrt{1+M^2}} d_c(p,q).$$
     We may assume $x_1\ge  x'_1$.   
  In this case,  the horizontal line segment $S=\{q*te_1|  t\in [0, x_1- x'_1]\}$
  has length $|x_1-x'_1|$  and 
    connects $q$ and $q'=q*(x_1-x'_1)e_1\in p*G$.  It follows that   $F(S)$ has length 
 $\le \sqrt{1+M^2}\cdot  |x_1-x'_1|$.  Hence 
$$d_c(F(q), F(q'))\le  \sqrt{1+M^2}\cdot  |x_1-x'_1|=\sqrt{1+M^2}\cdot d_c(q, p*G)\le  d_c(p,q)/10.$$
  On the other hand,  
  $d_c(F(p),  F(q'))=d_c(p, q')\ge (1-\frac{1}{10\sqrt{1+M^2}} ) d_c(p,q)$.  By the triangle inequality we have 
 $$d_c(F(p), F(q))\ge  d_c(F(p), F(q'))-d_c(F(q'), F(q))\ge  (1-\frac{1}{10\sqrt{1+L^2}}-\frac{1}{10} ) d_c(p,q).$$
  Hence $d_c(F(p), F(q))$ is bounded from below in terms of $d_c(p,q)$.  Since $F^{-1}_h=F_{-h}$,
  the same argument applied to $F_h^{-1}$   shows   that   
   $d_c(p, q)$ is  bounded from   below in terms of 
$d_c(F(p), F(q))$.  Hence $F$ is biLipschitz. 
 
   Conversely,   
  let $n\ge 3 $ and let $F:  F^n_\R\ra F^n _\R$ be a quasiconformal  map.  We have shown   that  after  
    composing $F$ with  graded automorphisms and left translations,     % $F(o)=o$  and 
   $F$ has  the following   form:
  $$F(x_1e_1*(\sum_{i=2}^{n+1}x_i e_i))=x_1e_1 *\sum_{i=2}^{n+1} \{(x_i+h_i(x_1)\}e_i.$$
   %  It remains  to 
   We next show that $h_2$  is Lipschitz.  % and that   $F=F_{h_2}$.
  %$h_{j}(x)=-\int_0^x  h_{j-1}(s)ds$ for $3\le j\le n+1$.
    Given $x_1, x'_1\in \R$, let $p=x_1e_1$, $q=x'_1e_1$.  Then $d(p,q)=|x_1-x'_1|$.  
   Since we have shown that $F$ is biLipschitz,   we have 
  $d(F(p), F(q))\le M\cdot d(p,q)=M\cdot |x_1-x'_1|$.    On the other hand,  
  the coefficient of $e_2$   in 
    $(-F(q))*F(p)$   %in the proof of Lemma \ref{f'} one easily sees that the coefficient of $e_2$ 
     is $h_2(x_1)-h_2(x'_1)$.  It follows that 
$$|h_2(x_1)-h_2(x'_1)|\le d(o,  (-F(q))*F(p))=d(F(q), F(p))\le   M\cdot |x_1-x'_1|.$$ 
  Hence $h_2$  is Lipschitz.  
  Now consider the quasiconformal map $F_0:=F^{-1}_{h_2}\circ F=F_{-h_2}\circ F$.  
   Its projection on the first layer is the identity map.   That is, if $\pi_1: F^n_\R\ra V_1$ denotes 
 the projection onto the first layer, and if $g\in F^n_\R$ is such that $\pi_1(g)=x_1e_1+x_2e_2$,  then $\pi_1(F_0(g))=x_1e_1+x_2e_2$.  It follows that the Pansu differential of $F_0$ is the identity isomorphism whenever it exists.  Now  Lemma \ref{pansu=id} implies $F_0$ is a left translation.  Hence the original map $F$ is a composition of left translations, graded isomorphisms and a map of the form $F_h$.

 % Furthermore, the proof of Lemma \ref{hi} shows that $-h_i(x_1)=h'_{i+1}(x_1)$. This implies 
   %$h_{i+1}(x_1)=h_{i+1}(0)-\int_0^{x_1} h_i(s)ds$.  Since we already normalized the map $F$ with $F(o)=o$,
  %we have $h_i(0)=0$. Hence    $h_{i+1}(x_1)=-\int_0^{x_1} h_i(s)ds$  and $F=F_h$.

The proof of Theorem \ref{thA} is now complete. 

\qed

\section{The complex Heisenberg  groups}\label{s4}

In this section we will   provide  evidence that quasiconformal maps on the complex Heisenberg groups are very special (Subsection \ref{QCcomplexh}).    % In fact we h
 For this purpose, we need to introduce differential forms associated with 
 $2$-step Carnot groups  (Subsection \ref{df}),    
   and    discuss  their relations with horizontal liftings (Subsection \ref{hlift}).

\subsection{Differential forms  associated with $2$-step Carnot groups}\label{df}

Here we 
introduce differential forms associated with 
 $2$-step Carnot groups.

Let $\mathcal N=V_1\oplus V_2$ be a 2-step Carnot group. Here we
identify the Lie group with its Lie algebra via the exponential map.
    The
Lie bracket restricted to the first layer  $V_1$  gives rise to a
 skew symmetric bilinear map
 $$\omega:  V_1\times V_1 \ra V_2,$$
  $$ \omega(X, Y)=[X,Y].$$
  We view  $\omega$ as a (constant) $V_2$-valued
 differential 2-form  on $V_1$.

We   next   define a $V_2$-valued differential 1-form $\alpha$ on $V_1$
 as follows.   For each $X\in V_1$, we need to define a linear map
$\alpha_X:  T_X V_1\ra V_2$. We identify  $T_X V_1$ with $V_1$.
  Let  $\alpha_X: V_1\ra V_2$ be given by
  $$ \alpha_X(Y)=[X,Y],\;\; \text{for}\;\; Y\in V_1. $$

It is convenient to write the differential forms $\alpha$ and
$\omega$ in coordinates. Fix a vector space basis $\{e_1, \cdots,
e_m\}$   for  $V_1$ and  a vector space basis $\{\eta_1, \cdots,
\eta_n\}$   for  $V_2$.  Then a point $X\in V_1$ can be written as
$X=x_1 e_1+ x_2 e_2+\cdots +x_m e_m$.
 For $X=x_1 e_1+ x_2 e_2+\cdots +x_m e_m$
   and $Y=y_1 e_1+ y_2 e_2+ \cdots+y_m e_m$,  one obtains:
     $$\omega(X, Y)=[X, Y]=\sum_{i,j} x_iy_j[e_i, e_j].  $$
      It follows that
        $$\omega=\sum_{i<j} [e_i, e_j]  dx_i\wedge dx_j$$
         and
         $$\alpha_X=\sum_{i, j} x_i\,[e_i, e_j]dx_j.$$
            One notices that   $d\alpha=2\omega$.

Here we work out the  differential forms associated to the complex Heisenberg group 
  $H^1_\C=\C^3$.  Let $X, Y, Z$ be the basis for the complex Lie algebra with bracket relation $[X,Y]=Z$.
    We choose  a basis for the real Lie albegra   $\{ e_1=X,  e_2=iX, e_3=Y, e_4=iY,  \eta_1=Z, \eta_2=iZ\}$.  
         The non-trivial bracket relations   for the real Lie algebra  
   are   $[e_1, e_3]=\eta_1$,   $[e_1, e_4]=\eta_2$,  $[e_2,  e_3]=\eta_2$  and $[e_2, e_4]=-\eta_1$.  
     A  point with  coordinates $(x_1, x_2, x_3, x_4)$ with respect to the real vector space  basis $\{e_1, e_2, e_3, e_4\}$
    has  coordinates $(w_1, w_2)$ with repect to the  complex vector space basis 
 $\{X,   Y\}$, where $w_1=x_1+ix_2, w_2=x_3+ix_4$.  
    Then $dw_1=dx_1+i dx_2$  and $dw_2=dx_3+i dx_4$.   
      Now 
  we obtain:  
  \b{align*}
\omega & =[e_1, e_3] dx_1\wedge dx_3+ [e_1,  e_4]  dx_1\wedge dx_4+[e_2,  e_3] dx_2\wedge dx_3
 +[e_2,  e_4]dx_2\wedge dx_4\\
&= (dx_1\wedge dx_3-dx_2\wedge dx_4)\eta_1+(dx_1\wedge dx_4 +dx_2\wedge dx_3) \eta_2\\
&=[ (dx_1\wedge dx_3-dx_2\wedge dx_4)+i (dx_1\wedge dx_4 +dx_2\wedge dx_3)] Z\\
&=(dw_1\wedge dw_2)Z.
\end{align*}
  Similarly,
   \b{align*}
\alpha  & =[e_1, e_3]x_1dx_3+[e_3, e_1]x_3dx_1 +[e_1, e_4]x_1dx_4+[e_4, e_1]x_4dx_1+\\
& +[e_2, e_3]x_2dx_3+[e_3, e_2]x_3dx_2 +[e_2, e_4]x_2dx_4+[e_4, e_2]x_4dx_2\\
&=(x_1dx_3-x_3dx_1+x_4dx_2-x_2dx_4)\eta_1
+(x_1dx_4-x_4dx_1+x_2dx_3-x_3dx_2)\eta_2 \\
&=[(x_1dx_3-x_3dx_1+x_4dx_2-x_2dx_4)
+i (x_1dx_4-x_4dx_1+x_2dx_3-x_3dx_2)]Z\\
&=(w_1dw_2-w_2dw_1)Z
\end{align*}

\subsection{Horizontal lifts in  $2$-step Carnot groups}\label{hlift}

Here we give a criteria for a closed curve in $V_1$   whose  horizontal
lifts  to $\mathcal N$  are  also  closed curves.

We  first  recall that a horizontal curve is completely determined by
its initial point and
  its first layer component.  This result is well-known.

Let  $c(t)=(c_1(t), c_2(t))\in \mathcal N=V_1\oplus V_2$ be an
absolutely continuous
   curve in  $\mathcal N$.  For each $t_0$, let
   $\gamma(t)=L_{-c(t_0)} c(t)$ be the translated curve. Notice that
   $\gamma(t_0)=0$.
 Using BCH formula, one finds that the  tangent vector of $\gamma$
 at $t=t_0$  is
  $$\big(c'_1(t_0), \;c'_2(t_0)-\frac{1}{2}[c_1(t_0), c'(t_0)]\big).$$
    It follows that the curve $c(t)$ is horizontal if and only
    if
     \b{equation}\label{e21}
     c'_2(t)=\frac{1}{2}[c_1(t), c'_1(t)]\;\; \text{for a.e. }\;
       t.
\end{equation}

%\subsection{Horizontal Lifts}\label{s3}

  Here is  a criteria for a closed curve in $V_1$   to have closed   horizontal
lifts  to $\mathcal N$.

\b{Le}\label{l1}
     Let $c_1: [0,1]\ra V_1$ be a closed  Lipschitz
         curve.  Then the following conditions are equivalent:\newline
  (1) The    horizontal lifts of $c_1$  to $\mathcal N$    are
closed curves;\newline
  (2)  $\int_{c_1}\alpha=0$;\newline
    (3)  $\int_D \omega=0$
    for any Lipschitz  2-disk  $D$ with boundary curve $c_1$.

\end{Le}

\b{proof}

$(1)\iff(2)$
    Let  $c(t)=(c_1(t), c_2(t))$ be  a  horizontal lift of
  $c_1$.   Then $c$ is closed if and only if $c_2(1)=c_2(0)$.   By (\ref{e21}),
   we have $c'_2(t)=\frac{1}{2} [c_1(t), c'_1(t)]$ for a.e.   $t\in
   [0,1]$.
   Write  $c_1(t)=\sum_i x_i(t) e_i$.  Then
     $c'_1(t)=\sum_i x'_i(t) e_i$  and
     $[c_1(t), c'_1(t)]=\sum_{i,j}x_i(t)x'_j(t)[e_i, e_j]$.
 Now the fundamental theorem
of calculus gives
   $$c_2(1)-c_2(0)=\int_0^1 c'_2(t) dt=\frac{1}{2} \int_0^1 [c_1(t), c'_1(t)] dt=
    \frac{1}{2} \int_0^1\sum_{i,j}x_i(t)x'_j(t)[e_i,
    e_j]\,dt=\frac{1}{2}\int_{c_1} \alpha .$$
      Hence (1) and (2)  are
    equivalent.

     (2) and (3) are equivalent due to  Stokes'    theorem  and the fact that $d\alpha=2\omega$.

\end{proof}

  \subsection{Quasiconformal maps on the complex Heisenberg groups}\label{QCcomplexh}

  Here 
we    provide  evidence that quasiconformal maps on the complex Heisenberg groups are very special.

Recall that the $n$-th complex Heisenberg group $H^n_\C$  is the simply connected Lie group whose Lie algebra   $\mathcal{H}^n_\C$  is a complex Lie algebra and   has a complex vector space basis
  $X_i, Y_i, Z$ ($1\le  i\le n$)  with the only non-trivial bracket relations 
  $[X_i, Y_i]=Z$,   $1\le i\le n$.   Of course, it has  more bracket relations as a real Lie albegra coming from the fact that the bracket is complex linear:   $[X_j, iY_j]=[iX_j,   Y_j]=i Z$  and $[iX_j,  iY_j]=-Z$.
  %one has bracket relations like $[X_i, iY_i]=i Z$ that come from the fact the bracket is complex linear.  
 The first layer $V_1$ of 
  $\mathcal{H}^n_\C$   is spanned by the $X_i, Y_i$, $1\le i\le n$ and has   complex  dimension $2n$. The second layer  $V_2$   is spanned by $Z$ and has complex  dimension $1$. 
    We identify   both  $H^n_\C$   and  its  Lie algebra   $\mathcal{H}^n_\C$   with $\C^{2n+1}$.  
  So $V_1=\C^{2n}\times \{0\}$ and $V_2=\{0\}\times \C$.    Let $\pi_1:  
{H}^n_\C=\C^{2n+1}   \ra \C^{2n}$ be the projection onto $V_1$.  
  A map $F: \C^{2n+1}   \ra \C^{2n+1} $  is called a lifting of a map $f: \C^{2n }  \ra \C^{2n}$
  if $F(\pi_1^{-1}(p))=\pi_1^{-1}(f(p))$  for all $p\in \C^{2n}$.   

  Let $\tau:  H^n_\C\ra H^n_\C$  be defined by  $\tau(w_1, \cdots, z)=(\bar{w_1}, \cdots, \bar{z})$.  
  It is easy to see that  $\tau$ is a graded isomorphism of $H^n_\C$.

Here is the first evidence for  Conjecture  \ref{conjectyre-ch}.

\b{Prop}\label{affine}
Let $F: H^n_\C\ra H^n_\C$   be a homeomorphism of the complex Heisenberg group. If $F$ is 
 both a quasiconformal map and a $C^2$ diffeomorhism, then 
  after possibly composing with $\tau$,   $F$ is the lifting of a  complex  affine map.
  Furthermore, $F$   or $F\circ \tau$ is  a complex  affine map.

\end{Prop}

\b{proof}   Write $F$ as $F(w_1, \cdots, w_{2n}, z)=(F_1, \cdots, F_{2n}, F_{2n+1})$, where $F_i=F_i(w_1, \cdots, w_{2n}, z)$ is the $i$-th component function of $F$.    
   Since $F$ is quasiconformal  and $H^n_\C$ is connected, $F$ is power quasisymmetric.  
  This implies that the $F_i$'s   have polynomial growth.  On the other hand,
by   the main result in \cite{RR},     $F$   or  $F\circ \tau$   is biholomorphic.      
    We shall assume $F$ is  biholomorphic.    
    By a Liuville type theorem, we  conclude that  the  $F_i$'s
     are actually polynomials.   By symmetry, the component functions of $F^{-1}$ are also polynomials. 
    It follows that the horizontal  Jacobian  $J_H(F, x)$  of $F$ is a polynomial.  The same is true for the horizontal Jacobian  $J_H(F^{-1}, y)$  of $F^{-1}$.  However, the horizontal Jacobian of $Id=F^{-1}\circ  F $ is
  $1=J_H(F, x)\cdot J_H(F^{-1}, F(x))$. 
   So the product of the polynomials    $J_H(F, x)$    and $ J_H(F^{-1}, F(x))$ is $1$.  
   This happens only when both polynomials are constants. 
   Therefore   $J_H(F, x)$ is a constant function.  
   This implies that $F$ is biLipschitz.    It follows that there is a constant $L>0$ such  that  
   $d(F(p), F(q))\le  L\cdot d(p, q)$ for all $p, q\in H^n_\C$.    Let $p=o$ and $q=(w_1, \cdots, w_{2n}, z)$.
  We see that  
 % $$|F_{2n+1}(w_1, \cdots, w_{2n}, z)-F_{2n+1}(0,\cdots, 0)|{\frac{1}{2}}\le L\cdot \big\{\sum_i |w_i|+ |z|^{\frac{1}{2}}\big\}$$   and  
    for each $1\le i\le 2n$, 
  $$|F_i(w_1, \cdots, w_{2n}, z)-F_i(0,\cdots, 0)|\le L\cdot \big\{\sum_i |w_i|+ |z|^{\frac{1}{2}}\big\}.$$
    Since $F_i$ is a   polynomial,    we conclude that  $F_i$ is independent of $z$ and is affine in $w_1, \cdots, w_{2n}$.  
  Hence $F$ is the lifting of a complex affine map  $f:=(F_1, \cdots, F_{2n}):    \C^{2n}\ra \C^{2n}$. 

Let $g$ be the linear part of $f$. Notice that  $g=dF(p)|_{V_1}$  for any point $p$.  In other words,  
   $g$ lifts to the  graded isomorphism   $dF(p)$    of  $H^n_\C$.  The translational part of $f$ of course 
 lifts  to a left translation   $L_q$  (for some $q\in H^n_\C$)   in $H^n_\C$.  Let $F'=L_q\circ dF(p)$.
   Notice that $dF(x)|_{V_1}=g=dF'(x)|_{V_1}$ for all $x\in  H^n_\C$.   It follows that 
 $dF(x)=dF'(x)$ for all $x\in  H^n_\C$.  By Lemma \ref{pansu=id},
  there is some $q'\in H^n_\C$  such that $F=L_{q'}\circ F'=L_{q'}\circ  L_q\circ dF(p)$,
     which is a complex affine map.

\end{proof}

Here  is   more  evidence for  Conjecture  \ref{conjectyre-ch}.

\b{Prop}\label{affine3}
 Let $F: H^1_\C\ra H^1_\C$  be a quasiconformal map.   Suppose $F$ is the lifting of a map  $f:
\C^2\ra \C^2$  
   of the form
 $f(w_1, w_2)=(w_1, w_2+g(w_1))$, where $g: \C\ra \C$ is a map.  
  Then there are constants $a, b\in \C$ such that   $g(w_1)=a w_1+b$.

\end{Prop}

\b{proof}   Notice that for each fixed $w_1\in \C$ we have $f(\{w_1\}\times \C)=\{w_1\}\times \C$.  
Since $F$ is a lifting of the map $f$, it follows that $F$ preserves each left coset of 
 $\{0\}\times \C^2$  in $H^1_\C$.  Now the arguments in \cite{SX}
    show that $F$ is biLipschitz.

Fix $w_2, z\in \C$ and  let  $c: [0,1]\ra \C\times \{w_2\}\times \{z\}\subset H^1_\C$
be a closed   $C^1$ curve. 
  Then $c$ is a horizontal curve  in  $H^1_\C$.     Being Lipschitz,     $F\circ c$ is   also a  closed    horizontal     curve in $H^1_\C$.
  The projection
  $\pi_1\circ F\circ c$  of $F\circ c$  under $\pi_1$ is a closed   Lipschitz  curve in $\C^2$ that admits a closed horizontal lift.   By Lemma   \ref{l1}    we have 
   $$\int_{\pi_1\circ F\circ c} w_1dw_2-w_2dw_1=0.$$
  Since $\int_\gamma w_1dw_2+w_2dw_1=0$   holds  for any closed curve $\gamma$, we have 
  $\int_{\pi_1\circ F\circ c} w_2dw_1=0.$
  Notice that $\pi_1\circ F\circ c(t)=(c(t),  w_2+g(c(t)))$.
  So  $$0=\int_{\pi_1\circ F\circ c} w_2dw_1=\int_c [w_2+g(w_1)]dw_1=\int_c g(w_1)dw_1.$$
   Hence $\int_c g(w_1)dw_1=0$ for any closed   $C^1$  curve in the complex plane.  
  Since $g$ is continuous (actually Lipschitz, see next paragraph),   
    Morera's theorem implies that $g(w_1)$ is holomorphic. 

We next show that $g$ is Lipschitz.     Let $w_1, w'_1\in \C$ be arbitrary.  
  Fix any $w_2, z\in \C$  and  let $p=(w_1, w_2, z)$,   $q=(w'_1, w_2,  z+\frac{1}{2}(w_1-w'_1)w_2)$.  
  Then $d(p, q)=|w_1-w'_1|$.  Notice
  $\pi_1\circ F(p)=f(w_1, w_2)=(w_1, w_2+g(w_1))$ and 
$\pi_1\circ F(q)=f(w'_1, w_2)=(w'_1, w_2+g(w'_1))$.
  It follows that $d(F(p), F(q))\ge |(w_2+g(w_1))-(w_2+g(w'_1))|=|g(w_1)-g(w'_1)|$.  
By the first paragraph,  
$F$ is $M$-biLipschitz  for some $M>0$. 
 Hence
 $$|g(w_1)-g(w'_1)|\le d(F(p), F(q))\le M\cdot d(p,q)=M\cdot |w_1-w'_1|.$$
     So  $g$ is  Lipschitz.   Since $g$ is also holomorphic,  it has to be affine.

\end{proof}

\section{Quasiconformal maps on the  complex model Filiform groups}\label{s3}

 In this Section we show that  quasiconformal maps on the higher complex 
  model Filiform groups are even more special than in the real case.   
   The proof is mostly similar to the real case. We will only indicate the difference in the proofs.

\subsection{Graded automorphisms of $\mathcal{F}^n_\C$}\label{gradedcomplex}

In this  Subsection we identify the graded  automorphisms of $\mathcal{F}^n_\C$.

%  We  first  show that graded isomorphisms
 %of the  complex model Filiform algebras  are either complex linear or complex   anti-linear.

Let  $F_\C^n$     be the $n$-step complex model Filiform group.   Recal that 
     its Lie algebra 
    $\mathcal{F}_\C^n$  
is a complex Lie algebra with basis $\{e_1, e_2, \cdots, e_{n+1}\}$
  and the only non-trivial bracket relations   are 
    $[e_1, e_j]=e_{j+1}$  for $2\le j\le n$.  % and all other brackets are zero. 
  Viewed as a real Lie algebra,    $\mathcal{F}_\C^n$  has the additional   bracket relations: 
%  In particular since the brackets are complex linear, we also have 
 $[ie_1, e_j]=ie_{j+1}=[e_1, ie_j]$,  
 $[ie_1, ie_j]=-e_{j+1}$.  % and $[e_j, ie_1]=i[e_1, e_1]=0$. 
   The    Lie algebra 
    $\mathcal{F}_\C^n$  decomposes as 
    $\mathcal{F}_\C^n=V_1\oplus V_2 \oplus\cdots \oplus V_n$, 
  where $V_1=\C e_1\oplus \C e_2$ is the first layer  and $V_j=\C e_{j+1}$,  $2\le j\le n$
  is the $j$-th layer.

The proof of the following  Lemma   is similar to that of Lemma \ref{l2.1}.

\b{Le}\label{l3.1}
Let $\mathcal{F}_\C^n$ be the $n$-step complex  model Filiform algebra.
  Assume $n\ge 3$.   Let $X=ae_1+be_2\in V_1$  be a  nonzero  element in the first layer. 
  Then   $\text{rank}(X)=2$  if $a=0$  and $\text{rank}(X)>2$  otherwise.

\end{Le}

  Hence,

\b{Le}\label{l3.2}
  Suppose $n\ge 3$.  Then  $h(\C e_2)=\C e_2$ for   every graded isomorphism $h: \mathcal{F}_\C^n
   \ra \mathcal{F}_\C^n$.

\end{Le}

A  graded isomorphism $h: \mathcal{F}_\C^n
   \ra \mathcal{F}_\C^n$  is in particular an isomorphism bewteen real vector spaces. In general, 
  $h$ needs not to be an isomorphism of the complex vector spaces.   Here is an example.   
  Define $\tau:  \mathcal{F}_\C^n
   \ra \mathcal{F}_\C^n$    by
  $$\tau(\sum z_j e_j)=\sum \bar{z_j} e_j.$$
  It is easy to check that $\tau$ is a graded isomorphism of $\mathcal{F}_\C^n$.
  Notice that $\tau$ is complex anti-linear, not complex linear.

\b{Le}\label{complel}
      Let   $h: \mathcal{F}_\C^n
   \ra \mathcal{F}_\C^n$  be a 
 graded isomorphism.
  Then $h$ is either complex linear or complex anti-linear.

\end{Le}

\b{proof}  The map $h$ 
 satisfies   $[hv, hw]=h[v,w]$  for all $v, w\in  V_1$.    In particular, we have 
  $[h(ie_1), h(i e_2)]=h([ie_1, ie_2])=-h([e_1, e_2])$,
  $[h(ie_1), h(e_2)]=h([ie_1, e_2])=h([e_1, ie_2])=[h(e_1), h(ie_2)]$  and 
  $[h(e_1), h(ie_1)]=h([e_1, ie_1])=0$.

 By    Lemma \ref{l3.2}   we know that $A(\C e_2)=\C e_2$.   
  Hence there are $0\not=w_1,    w_2\in \C$ such that 
 $$h(e_2)=w_1  e_2$$
  $$ h(ie_2)=w_2e_2.$$
    There are also constants $a, b, c,d\in   \C$ 
  such that 
$$h(e_1)=a e_1+b  e_2$$
  $$h(ie_1)=c  e_1+ d  e_2.$$

The equation 
 $[h(ie_1), h(e_2)]=[h(e_1), h(ie_2)]$  yields   $cw_1=a w_2$.  
  Similarly,   $[h(ie_1), h(ie_2)]=-[h(e_1), h(e_2)]$ yields
  $cw_2=-aw_1$  and $[h(e_1), h(ie_1)]=0$ yields  $ad-bc=0$.  
  It follows that $c(w_1^2+w_2^2)=0$.  
    Assume   $c=0$, then  $ad=0$.  
 Since $h$ is an isomorphism,   $d\not=0$.   So $a=0$,    which implies that $h$ maps $V_1$ into 
 $\C e_2$, contradicting the fact that $h$ is an isomorphism.    Hence  $c\not=0$.
  It follows that $w_2=iw_1$  or $w_2=-iw_1$.

If $w_2=iw_1$,  then $c=ia$  and $d=ib$. In this case,  $h|_{V_1}$  is complex linear. 
   Since the bracket is  also complex linear,  an induction argument shows 
 $h|_{V_j}$ is complex linear for all $j$.  
 If $w_2=-iw_1$,  then 
  $c=-ia$  and $d=-ib$.  In this case,  $h|_{V_1}$  is complex anti-linear.   A similar
  argument also shows that 
 $h|_{V_j}$ is complex   anti-linear for all $j$.  

\end{proof}

\b{Le}\label{graded-complex}
  A real linear map  $h: \mathcal{F}_\C^n
   \ra \mathcal{F}_\C^n$   is a graded isomorphism if  and only if $h$   has one of the following two forms:\newline
  (1)  $h=h_{a_1, a_2, b}$  for some $a_1, a_2\in \C\backslash\{0\}$ and $b\in \C$; in this case, $h$ is complex linear;\newline
  (2)  $h=\tau\circ h_{a_1, a_2, b}$  for some $a_1, a_2\in \C\backslash\{0\}$ and $b\in \C$; in this case, $h$ is complex   anti-linear.

\end{Le}

\b{proof}
 One direction is clear. So we start with a 
graded isomorphism  
$h: \mathcal{F}_\C^n
   \ra \mathcal{F}_\C^n$.   By Lemma \ref{complel}  $h$ is complex linear or complex  anti-linear.  
  If $h$ is complex linear, then the proof of Lemma \ref{graded-real}  shows that $h=h_{a_1, a_2, b}$  for some $a_1, a_2\in \C\backslash\{0\}$ and $b\in \C$.  If  $h$ is complex anti-linear, then $\tau\circ h$ is 
  a graded isomorphism and is complex linear. Hence 
  $\tau\circ h=h_{a_1, a_2, b}$  for some $a_1, a_2\in \C\backslash\{0\}$ and $b\in \C$. 
  It follows   that  $h=\tau\circ h_{a_1, a_2, b}$.

\end{proof}

\subsection{Proof of Theorem \ref{thB}}\label{specialcomplex}

In this  Subsection we show that each quasiconformal   map of  $F^n_\C$  ($n\ge 3$) is a  composition of 
  left translations  and graded automorphisms.    Hence by Lemma \ref{graded-complex},
   after possibly  taking    complex  conjugation   (that is, composing with $\tau$),  each quasiconformal map 
of  $F^n_\C$  ($n\ge 3$)  is  biholomorphic.

  Let  $n\ge 3$  and $F: {F}_\C^n  \ra  {F}_\C^n$   be a quasiconformal map.  Then $F$ is a quasisymmetric map
  and   sends 
  left cosets of  $\C e_2$ to  left cosets of $\C  e_2$.  
  For   a.e. left coset
  $L$  of $\C e_2$,  the map $F$ is Pansu differentiable a.e. on $L$.    
    For any point $x\in F_\C^n$  where $F$ is Pansu differentiable,   
  let $a_j(x),  b(x)\in \C$,  $1\le j\le n+1$,   be such that 
$dF(x)(e_1)=a_1(x) e_1+ b(x) e_2$,     $dF(x)(e_j)=a_j(x) e_j$ ($2\le j\le n+1$).  
  Then  $a_j(x)=(a_1(x))^{j-2}a_2(x)$ for $2\le j\le n+1$. 
   % Then   $dF(x)(e_3)=a(x) d(x) e_3$.   We denote $a_2(x)=a(x) d(x)$.  
  If $dF(x)$  is complex  linear, then $dF(x)(ze_j)=z\cdot dF(x)(e_j)$  for any $z\in \C$.
 If $dF(x)$  is complex  anti-linear, then $dF(x)(ze_j)=\bar{z}\cdot  dF(x)(e_j)$  for any $z\in \C$.

\b{Le}\label{secondcomplex}
Let $L$ be a left coset of $\C e_2$ where $F$ is Pansu  differentiable a.e.
  Then there is a constant $0\not=a_L\in \C$  such that $a_3(x)=a_L$  for a.e. $x\in L$. 
 Furthermore,  if $x, y\in L$ are two points where $F$ is Pansu-differentiable,
  then  either both  $dF(x)$ and $dF(y)$ are complex linear or both are complex  anti-linear.

\end{Le}

\b{proof} 
 The proof of the first statement  is the same as in the proof of  Lemma \ref{second}
      with $\R$ replaced by   $\C$.  
  For the second   statement,  if      $dF(x)$  is  complex linear  and $dF(y)$ is complex anti-linear,
  then use  the left coset $ir^2e_3*\C e_2$ instead of $r^2e_3*\C e_2$  and repeat the argument to get a contradiction.   For this,  one uses the fact that $dF(x)(ir^2 e_3)= ir^2 a_3(x) e_3$  and 
  $dF(y)(ir^2 e_3)= -ir^2 a_3(y) e_3=-ir^2 a_3(x) e_3$.

\end{proof}

Now  statements  similar to  
Lemma   \ref{new1}   through  Lemma \ref{e2}  all hold and the proofs are   the same.  
  In particular,  $F$ is biLipschitz, and if we 
 denote   by $G=\C e_2\oplus \cdots \oplus \C e_{n+1}$, then   $F$ sends left cosets of $G$ to left cosets of $G$.  Hence  there is a homeomorphism $f: \C\ra \C$ such that $F(x_1e_1*G)=f(x_1)e_1*G$.  
For a.e. left coset $L$  of  $\C e_2$,    there exists a constant $0\not=a_{2, L}\in \C$ such that 
 $a_2(x)=a_{2,L}$  for a.e. $x\in L$.     Furthermore, for any $g\in L$,    either 
$F(g* te_2)=F(g)* a_{2,L} te_2$    for all  $t\in \C$  or  
 $F(g* te_2)=F(g)* a_{2, L}\bar{t}e_2$  for all $t\in \C$.

%Now Lemma  \ref{secondcomplex} implies that 
  %for each left coset $L:=g\cdot \C e_2$ as in the lemma,  
   %there is a nonzero constant $a_L\in \C$ such   that  
    %$F(g* te_2)=F(g)* a_Lte_2$   for all $t\in \C$   or   
  %$F(g* te_2)=F(g)* a_L\bar{t}e_2$  for all $t\in \C$. 

\b{Le}\label{same2}
 There is a constant $0\not=a_2\in \C$ such that 
 either $F(g* te_2)=F(g)* a_2te_2$  holds for all $g\in F_\C^n$ or 
$F(g* te_2)=F(g)* a_2\bar{t}e_2$   holds for 
  all $g\in F_\C^n$.  

\end{Le}

\b{proof}
Let $L,  L'$ be two  left cosets of $\C e_2$ where $F$ is Pansu  differentiable a.e.
We run the argument in  the proof of Lemma \ref{same}      for $t\in \R$  to show $a_{2,L}=a_{2, L'}$.  
   If  $F(g* te_2)=F(g)* a_{2,L}te_2$  holds on  $L$  and  $F(g* te_2)=F(g)* a_{2, L'}\bar{t}e_2$  holds on      $L'$,  then 
   the same arguemnt using $t\in i\R$
    yields   a contradiction. 
   Since $F$ is a homeomorphism  and   for  a.e. left coset $L$  of $\C e_2$,   the map $F$ is Pansu differentiable a.e. on $L$,   the lemma follows by continuity.

\end{proof}

After composing $F$ with  a graded isomorphism, we may assume that   $a_2=1$  and 
$F(g* te_2)=F(g)* te_2$  holds for all $g\in F_\C^n$.  So the Pansu differential is complex linear  whenever it exists. 

When the Pansu differential is complex linear,  the derivatives that appear  in  the proofs in  Section \ref{s2}  
  can be taken  to be  complex derivatives:\newline
 (1) proof of Lemma \ref{e2},  $a_1(x)=f'(x_1)$;\newline
  (2)  proof of Lemma \ref{add1},  $\frac{\partial f_j}{\partial x_j}(p)=(f'(x_1))^{j-2}$; here both are complex derivatives;  \newline
 (3)  proof of Lemma \ref{hi},  $-h_n(x_1, x_{n+1})=h'_{n+1}(x_1)$.

     % Let $p=\sum_i x_ie_i$ be a point where $F$ is Pansu differentiable. Then proposition ?? and Lemma ? imply   that the Pansu differential at $p$ has the form:  
  %$dF(p)(e_1)=b(x_1) e_1 +c(p) e_2$,    $dF(p)(e_2)=a e_2$,   $dF(p)(e_i)=b(x_1)^{i-2}a e_i$  for $3\le i\le n$.     After comp[osing $F$ with a 
 %graded isomorphism we  may assume $a=1$.  

 Notice that the left cosets  of $G$ are  isometric to 
 $\C^n$ with the metric 
 $D((z_i), (w_i))=\sum_i |z_i-w_i|^{\frac{1}{i}}$.  
   By the proof in Section  15 of    \cite{T}     each quasisymmetric map $h: (\C^n, D)\ra (\C^n, D)$ preserve the foliation
    consisting of    affine subspaces parallel to $\C^i\times \{0\}$  for each  $1\le i \le n-1$.  
  An analogue of Lemma  \ref{add1} holds, and so  for a.e.   $x_1\in \C$,  $F$ has the following form:
 $$F(x_1e_1*(\sum_{i=2}^{n+1}x_i e_i))=f(x_1)e_1 *\sum_{i=2}^{n+1} \{(f'(x_1))^{i-2}x_i+h_i(x_1,   x_{i+1},  \cdots, x_{n+1})\}e_i.$$

  We next show:

\b{Le}\label{fcomplexa}
The function $f$ is complex affine.
\end{Le}

\b{proof}  Recall that   $f: \C\ra \C$
 is a homeomorphism 
 such that $F(x_1e_1*G)=f(x_1)e_1*G$.  
  %We note that $F$ sends left cosets of $G$ to left cosets of $G$.  
  This implies that at any point $p\in F^n_\C$ where $F$ is Pansu differentiable, 
    $a_1(p)=f'(x_1)$.  Here $x_1$ is the  coefficient of $e_1$ in the expression for $p$ and $f'(x_1)$ is the complex derivative.   Hence at a.e. $x_1 \in \C$,  $f$ has nonzero complex derivative. In particular,  
 $f: \C\ra \C$ is a $1$-quasiconformal map.  
  %  By a result of Tukia-Vaisala \cite{TV},  
It follows that  $f$ is a similarity.  Since $f$ has complex derivative, 
    the linear part of $f$ can not be a reflection and so $f$ must be a complex affine map.

 % In this case, the arguments in \cite{SX}   imply  that 
   %$f$  is biLipschitz.    On the other hand, the proof of Lemma \ref{e2}
 %shows at a.e. $x_1\in \C$,  $f$ has complex derivative $f'(x_1)=a_1(x)$

\end{proof}

After composing $F$ with a graded isomorphism  we may assume  $f(x_1)=x_1$. So  $F$ has the following form:
$$F(x_1e_1*(\sum_{i=2}^{n+1}x_i e_i))=x_1e_1 *\sum_{i=2}^{n+1} \{(x_i+h_i(x_1,   x_{i+1}, \cdots, x_{n+1})\}e_i.$$
    Now  the proof of Lemma \ref{hi}   shows that $h_i$ is a function of $x_1$ only.  
We shall show that $h_i$ is a holomorphic function of $x_1$.  

Let  $H_3=\C e_4\oplus \cdots \oplus \C e_{n+1}$.  Then $H_3$ is a  closed  normal  subgroup of $F_\C^n$
  and $F_\C^n/{H_3}$ is isomorphic to    $H^1_\C$ (the first complex Heisenberg group).  
    It is easy to see  from the expression of $F$ that $F$ maps left cosets of $H_3$ to left cosets of 
  $H_3$.  Hence $F$ induces a map $\overline{F}:  F_\C^n/{H_3}=H_\C^1\ra  H_\C^1= F_\C^n/{H_3}$
  and $\overline{F}$ admits the following expression:
  $$\overline{F}(x_1e_1*(x_2e_2+x_3e_3))=x_1e_1*[(x_2+h_2(x_1))e_2+(x_3+h_3(x_1))e_3].$$
     %   The arguments in \cite{SX}   show that $F$ is biLipschitz.   
   It follows that $\overline F$ is the lifting of  the map
 $f: \C^2\ra \C^2$,   $f(x_1, x_2)=(x_1, x_2+ h_2(x_1))$.   
  Here we identified the first layer $V_1$ of  $H^1_\C$   
with $\C^2$ via  $x_1e_1+x_2e_2\ra (x_1, x_2)$.  

Let $\pi:  F_\C^n\ra F_\C^n/{H_3}=H^1_\C$ be the quotient map.  
Since $H_3$ is normal in $F_\C^n$,  and the quotient group 
$F_\C^n/{H_3}$  is also Carnot,  it is not  hard to check that 
 for any $p, q\in F_\C^n/{H_3}=H^1_\C$  and any  $x\in \pi^{-1}(p)$   one has 
$d_c(\pi^{-1}(p), \pi^{-1}(q))=d_c(x, \pi^{-1}(q))=d_c(p, q)$.   
  Since $F$ is quasisymmetric,  it now follows from the following Lemma of Tyson that 
$\overline{F}:  F_\C^n/{H_3}=H_\C^1\ra  H_\C^1= F_\C^n/{H_3}$  is also quasisymmetric.

\begin{Le}\label{tyson2} \e{(\cite[Lemma 15.9]{T})}
Let $g: X_1\ra X_2$  be an $\eta$-quasisymmetry    and $A,B,
C\subset X_1$.  If  $d(A,B)\le t\, d(A, C)$ for some $t\ge 0$, then
there is some $a\in A$ such that
\[
   d(g(A), g(B))\le \eta(t) d(g(a), g(C)).
\]
\end{Le}

 Now we apply  Proposition 
  \ref{affine3}   to  $\overline F$   to conclude  that $h_2(x_1)$ 
  is a  complex  affine  function of $x_1$.  
  So   there   are  constants   $a, b\in \C$ such that $h_2(x_1)=a x_1+b$.  After composing $F$ with the map   $F_{-h_2}$ (see Introduction)  we now may assume $h_2(x_1)=0$.    It follows that   the Pansu differential of $F$ is a.e. the identity isomorphism.   
 Now     Lemma  \ref{pansu=id}  implies that   it is  a   left translation.  Hence 
every quasiconformal map of  
 $F_\C^n$   is a finite composition of left translations and graded isomorphisms.   This finishes the proof of 
Theorem \ref{thB}.  

\qed

Notice that left translations in 
  $F_\C^n$   are polynomial  maps   with  polynomial inverse (this follows from the BCH formula).  By Lemma \ref{graded-complex},
 each graded isomorphism is complex linear after possibly composing with $\tau$ (taking complex conjugation).   Hence     Theorem \ref{thB}    implies that for $n\ge 3$,  every quasiconformal map
 of   $F_\C^n$   is a polynomial map with polynomial inverse after possibly composing with $\tau$.

 \addcontentsline{toc}{subsection}{References}

\noindent Xiangdong Xie: Dept. of Mathematics  and Statistics,  Bowling Green State  University,
     Bowling Green,   OH   43403,    U.S.A.\hskip .4cm E-mail:
xiex@bgsu.edu

\end{document}